\documentclass[reqno, 11pt]{amsart}
\usepackage{amssymb,bbm}
\usepackage{amsfonts}
\usepackage{amsmath, amscd}
\usepackage{graphicx}
\usepackage{color}
\usepackage{enumerate}
\usepackage{soul}
\usepackage{dsfont}
\usepackage{hyperref}
\hypersetup{
   colorlinks=true,
   linkcolor=blue,
   citecolor=blue,
   urlcolor=blue,
}

\usepackage{tikz}
\usetikzlibrary{automata}
\usetikzlibrary{arrows.meta,decorations.pathmorphing,decorations.pathreplacing,calc,positioning}

\makeatletter
\renewcommand\tocchapter[3]{%
  \indentlabel{\@ifnotempty{#2}{\ignorespaces#2.\quad}}#3%
}
\newcommand\@dotsep{4.5}
\def\@tocline#1#2#3#4#5#6#7{\relax
  \ifnum #1>\c@tocdepth
  \else
    \par \addpenalty\@secpenalty\addvspace{#2}%
    \begingroup \hyphenpenalty\@M
    \@ifempty{#4}{%
      \@tempdima\csname r@tocindent\number#1\endcsname\relax
    }{%
      \@tempdima#4\relax
    }%
    \parindent\z@ \leftskip#3\relax \advance\leftskip\@tempdima\relax
    \rightskip\@pnumwidth plus1em \parfillskip-\@pnumwidth
    #5\leavevmode\hskip-\@tempdima{#6}\nobreak
    \leaders\hbox{$\m@th\mkern \@dotsep mu\hbox{.}\mkern \@dotsep mu$}\hfill
    \nobreak
    \hbox to\@pnumwidth{\@tocpagenum{#7}}\par
    \nobreak
    \endgroup
  \fi}
\makeatother
\AtBeginDocument{%
\makeatletter
\expandafter\renewcommand\csname r@tocindent0\endcsname{0pt}
\makeatother
}
\def\l@subsection{\@tocline{2}{0pt}{2.5pc}{5pc}{}}

\newtheorem{theorem}{Theorem}[section]
\newtheorem{proposition}[theorem]{Proposition}
\newtheorem{lemma}[theorem]{Lemma}
\newtheorem{corollary}[theorem]{Corollary}
\newtheorem{definition}[theorem]{Definition}
\newtheorem{remark}{Remark}[section]

\newcommand{\Ind}[1]{\mathbbm{1}_{#1}}

\makeatletter
\@addtoreset{equation}{section}
\makeatother

\renewcommand\thetable{\thesection.\@arabic\c@table}

\title[The Contact Process with Temporary Immunity on $\mathbb{Z}$]{The Contact Process with Temporary Immunity on $\mathbb{Z}$}

\author[T. S. Mountford, R. Santos, C. Terra, M. E. Vares]{Thomas Mountford\textsuperscript{1}, Rafael Santos\textsuperscript{2}, C\'elio Terra\textsuperscript{2}, Maria Eul\'alia Vares\textsuperscript{2}} 

\address{\textsuperscript{1} Institute of Mathematics, \'Ecole Polytechnique F\'ed\'erale de Lausanne (EPFL), Lausanne, Switzerland} 
\address{\textsuperscript{2} Instituto de Matem\'atica, Universidade Federal do Rio de Janeiro (UFRJ), Rio de Janeiro, Brazil}
\email{thomas.mountford@epfl.ch}
\email{rafaels@dme.ufrj.br}
\email{eulalia@im.ufrj.br}
\email{caugusto.terra@gmail.com}

\begin{document}

\begin{abstract}
	We develop a regeneration scheme for the one-dimensional contact process with temporary immunity, also known as the spatial SIRS process, a non-attractive interacting particle system. The key ingredient is a one-sided space-time shielding mechanism, constructed from infinite infection paths of an auxiliary process. These paths separate the evolution of the infection front from the configuration and graphical history left behind, allowing us to construct regeneration times without relying on attractiveness. The same shielding mechanism also provides coupling arguments for processes started from different initial configurations. Using this approach, we establish a law of large numbers and a functional central limit theorem for the right edge, existence and uniqueness of an invariant law for the process seen from the edge, and complete convergence. We also prove survival from every nonempty finite infected set for sufficiently large immunity-loss rates and derive an explicit extinction criterion.
\end{abstract}

\keywords{interacting particle systems; CPTI; SIRS process; regeneration; right edge; functional central limit theorem; invariant measure}
\maketitle

\section{Introduction}
\label{sec:introduction}

Regeneration techniques provide a powerful approach to
the study of stochastic processes by identifying times
at which the relevant future evolution can be separated
from the preceding history. In the study of infection fronts, such methods have been developed through break-point constructions and coupling arguments for the classical contact process and related models; see, for instance,~\cite{CoxDurrettSchinazi1991,GalvesPresutti1987,Kuczek1989,MountfordSweet2000,Tzioufas2011}. These constructions exploit graphical properties that allow the future evolution of the front to be separated from its past.

The contact process with temporary immunity (CPTI),
also known as the spatial SIRS process, provides a
natural setting in which constructing such a separation
is particularly challenging. In this model, infected
individuals transmit infection to susceptible neighbors
and, upon recovery, enter a temporary immune state
before becoming susceptible again. Thus, previous
infections leave behind immune sites that can obstruct
subsequent transmission. Although the CPTI is Markovian, the legality of an infection path in its graphical construction depends on immunity intervals generated by earlier infections. Moreover, the model is non-attractive, so the monotone couplings underlying classical arguments are no longer available. The difficulty is therefore to recover a suitable separation from the past without relying on attractiveness.

In this paper, we develop a geometric alternative to monotone coupling based on a one-sided space-time shielding mechanism. Suitable infinite infection paths of an auxiliary process act as barriers that isolate the evolution on one side from modifications of the configuration and graphical history on the other. This allows us to replace the graphical history behind a shielding path without affecting the subsequent evolution of the infection front. By exploiting this property, we construct independent regenerative cycles for the front, even though the full process does not regenerate.

\subsection{Model and historical background}
\label{subsec:model}

The CPTI is a continuous-time Markov process on
$\{-1,0,1\}^{\mathbb Z}$, depending on two parameters
$\lambda,\mu>0$. The states $-1$, $0$, and $1$ represent
immune, susceptible, and infected individuals, respectively.
For a configuration $\xi\in\{-1,0,1\}^{\mathbb Z}$, let
\[
n_\xi(x):=\#\{y\in\mathbb Z:|y-x|=1,\ \xi(y)=1\}
\]
denote the number of infected neighbors of $x$. The transition
rates at site $x$ are
\[
0\to 1 \quad \text{at rate } \lambda n_\xi(x),\qquad
1\to -1 \quad \text{at rate } 1,\qquad
-1\to 0 \quad \text{at rate } \mu.
\]
Thus, $\lambda$ is the infection parameter and $\mu$ is the
rate at which immunity is lost. In the limiting case
$\mu=\infty$, the model reduces to the classical contact
process, while $\mu=0$ corresponds to permanent immunity,
yielding a spatial epidemic model of SIR type.

The presence of temporary immunity introduces a competition
between the propagation of infection and the creation of
immune sites. In particular, increasing the infection
parameter may facilitate the spread of infection while also
producing more immune sites that can obstruct its subsequent
propagation. As a consequence, the CPTI lacks the usual
attractiveness property, and monotonicity of its survival
probability with respect to the parameters is not known.

The CPTI provides a natural setting for studying the
interplay between spatial growth and local memory.
Infection promotes its own propagation through
transmission, but also creates temporary obstacles
that may hinder subsequent infections. Understanding
how these competing effects influence survival and
large-scale growth is a central question for the model.
Unlike contact processes in independently evolving
random environments, such as those studied in~\cite{Broman2007,LinkerRemenik2020,Remenik2008},
the immune environment is generated by the infection
itself. This also distinguishes the CPTI from the
renewal and dynamic-edge models considered in~\cite{HilarioUngarettiValesinVares2022}, where the
environmental restrictions evolve independently of
the infection. The resulting feedback between
infection and immunity introduces a history dependence
that is absent from these models.

The CPTI was introduced in~\cite{DurrettNeuhauser1991},
where the existence of a nontrivial invariant measure was
established in a suitable parameter regime. Subsequently,~\cite[Theorem~1.2]{AndjelSchinazi1996} extended the coexistence results to every dimension. In the same paper, they also established a complete convergence theorem in one dimension for a related but different model, in which the immunity-loss rate is proportional to the number of healthy nearest neighbors~\cite[Theorem~1.1]{AndjelSchinazi1996}.

By comparison with the classical contact process, extinction
from finite initial infections holds whenever the infection
rate is at most the critical value of the classical contact
process, independently of the immunity-loss rate.
Moreover,~\cite{vanDenBergGrimmettSchinazi1998} showed
that, for every fixed infection rate, the one-dimensional
CPTI dies out from finite initial infections when the
immunity-loss rate is sufficiently small.

A fundamental question, however, remained open in one
dimension: does the CPTI survive from a finite initial
infection for some choice of parameters? The existence
of a nontrivial invariant measure does not by itself
imply such survival in the absence of attractiveness.
In dimensions $d\geq 2$, survival from finite infections
was established for sufficiently large infection rates
in~\cite{vanDenBergGrimmettSchinazi1998}.
More recently,~\cite{Heeszel2026} showed that weak
survival implies strong survival for the one-dimensional
CPTI, while~\cite{HeeszelWascher2026} obtained survival
and extinction criteria in higher dimensions by adapting
a block construction originally developed for the contact
process with avoidance.

Numerical studies of the CPTI on lattices include
simulations and mean-field approximations
\cite{JooLebowitz2004}. On finite graphs, extinction
times have been studied for several graph classes
\cite{FriedrichEtAl2024,LamNguyenYang2024}.

Beyond the existence of nontrivial invariant measures,
the long-time behavior of the one-dimensional CPTI
has remained largely unexplored. Questions concerning
the propagation speed of the infection front, its
fluctuations, and convergence to equilibrium require
methods capable of controlling the dependence on past
infections. In the absence of attractiveness, standard
coupling and regeneration arguments cannot be applied
directly. The shielding mechanism developed here
overcomes this difficulty by isolating the evolution
relevant to the infection front from the graphical
history behind it.

\subsection{Main results}
\label{subsec:main-results}

For an initial configuration
$\xi\in\{-1,0,1\}^{\mathbb Z}$, let
$\mathbb P^\xi_{\lambda,\mu}$ denote the law of the CPTI
with parameters $(\lambda,\mu)$ started from $\xi$, and let
$\mathbb E^\xi_{\lambda,\mu}$ denote the corresponding
expectation. Throughout Sections~\ref{sec:block-renormalization-survival}-\ref{sec:shielding-asymptotic-laws} and Appendix~\ref{app:moment-estimates}, we fix $\lambda$ and suppress its dependence from the notation. When the initial configuration is clear from the context, we also omit it from the notation.

For $A\subset\mathbb Z$, let $\xi^A$ denote the process
started with the sites in $A$ infected and all other sites
susceptible, and write
\[
I_t^A:=\{x\in\mathbb Z:\xi_t^A(x)=1\}.
\]
Define the extinction time
\[
\tau^A:=\inf\{t\geq0:I_t^A=\varnothing\}.
\]
We say that the process survives from $A$ if
$\tau^A=\infty$, and set
\[
\theta(\lambda,\mu)
:=
\mathbb P_{\lambda,\mu}
\bigl(\tau^{\{0\}}=\infty\bigr).
\]

Let $\lambda_c(\infty)$ denote the critical infection rate
of the classical contact process. For fixed $\mu>0$, define
\[
\lambda_c(\mu)
:=
\inf\{\lambda>0:\theta(\lambda,\mu)>0\},
\]
with the convention $\inf\varnothing=\infty$.

Similarly, for fixed $\lambda>\lambda_c(\infty)$, define
\[
\mu_c(\lambda)
:=
\inf\{\mu>0:\theta(\lambda,\mu)>0\}
\]
and
\[
\overline\mu_c(\lambda)
:=
\sup\{\mu>0:\theta(\lambda,\mu)=0\}.
\]
Since monotonicity of $\theta(\lambda,\mu)$ in the
parameters is not known, these quantities are not
asserted to define a sharp phase transition curve.

Our first result gives an explicit sufficient condition
for extinction from finite initial infections, based on
the one-dimensional geometry of infection genealogies.

\begin{theorem}[Extinction criterion]
	\label{teo_dies}
	Let $\lambda,\mu>0$. If
	\begin{equation}
		\label{eq:extinction-criterion}
		\frac{\mu\lambda(1+2\lambda)}
		{(1+\mu)(1+\lambda)}<1,
	\end{equation}
	then, for every finite $A\subset\mathbb Z$,
	\[
	\mathbb P_{\lambda,\mu}(\tau^A<\infty)=1.
	\]
\end{theorem}

Our next result settles the question of survival from
finite initial infections in one dimension provided the immunity-loss rate is sufficiently large.

\begin{theorem}[Survival for large immunity-loss rate]
	\label{teo_survive}
	For every $\lambda>\lambda_c(\infty)$, there exists
	$\mu_0=\mu_0(\lambda)<\infty$ such that, for every
	$\mu\geq\mu_0$ and every nonempty finite
	$A\subset\mathbb Z$,
	\[
	\mathbb P_{\lambda,\mu}(\tau^A=\infty)>0.
	\]
	In particular, for every $\lambda>\lambda_c(\infty)$,
	\[
	\mu_c(\lambda)<\infty.
	\]
\end{theorem}

Together, these results extend the known one-dimensional
phase diagram in two directions: Theorem~\ref{teo_dies}
provides an explicit extinction criterion, while
Theorem~\ref{teo_survive} establishes a regime of
survival from finite initial infections.

Our main technical contribution is a one-sided space-time shielding mechanism that replaces key uses of attractiveness in classical coupling and regeneration arguments. The shielding paths allow us to modify the graphical history on one side of a space-time barrier without affecting the evolution on the other. This separation makes it possible to construct independent regenerative cycles for the infection front without regenerating the full configuration. The same mechanism also yields couplings between processes started from different initial configurations, leading to convergence to equilibrium.

Set
\[
\mathcal X^{\mathrm{e}}
:=
\left\{
\xi\in\{-1,0,1\}^{\mathbb Z}:
\xi(0)=1,\
\xi(x)\neq1\ \text{for all }x>0,\
\#\{x<0:\xi(x)=1\}=\infty
\right\}.
\]  
Thus, for $\xi\in\mathcal X^{\mathrm{e}}$, the origin
is the rightmost infected site and there are infinitely
many infected sites to its left; susceptible and immune
sites are otherwise allowed arbitrarily.

For $\xi\in\mathcal X^{\mathrm{e}}$, under
$\mathbb P^\xi_{\lambda,\mu}$, let
\[
r_t:=\sup\{x\in\mathbb Z:\xi_t(x)=1\}
\]
denote the position of the rightmost infected site
at time $t$.

Our first application of the regenerative structure
is a law of large numbers for the position of the
right edge.

\begin{theorem}[Law of large numbers for the right edge]
	\label{teo_lln}
	Fix $\lambda>\lambda_c(\infty)$. Then, for every
	$\mu\geq\mu_0(\lambda)$, there exists a deterministic
	constant $v=v(\lambda,\mu)>0$ such that, for every
	$\xi\in\mathcal X^{\mathrm{e}}$,
	\begin{equation}
		\label{eq:right-edge-lln}
		\frac{r_t}{t}\longrightarrow v
		\qquad
		\mathbb P^\xi_{\lambda,\mu}\text{-almost surely
			as }t\to\infty.
	\end{equation}
	In particular, the asymptotic speed does not
	depend on the initial configuration $\xi$.
	
	If $\xi$ is the configuration given by
	\[
	\xi(x)=
	\begin{cases}
		1, & x\leq0,\\
		0, & x>0,
	\end{cases}
	\]
	then the convergence in \eqref{eq:right-edge-lln}
	also holds in $L^1(\mathbb P^\xi_{\lambda,\mu})$.
\end{theorem}

As the immunity-loss rate tends to infinity,
the asymptotic velocity approaches that of the
classical contact process.

\begin{proposition}[Classical limit of the speed]
	\label{prop:speed-limit-mu-infty}
	Fix $\lambda>\lambda_c(\infty)$. Then
	\[
	\lim_{\mu\to\infty}v(\lambda,\mu)
	=
	v_{\rm CP}(\lambda),
	\]
	where $v_{\rm CP}(\lambda)$ is the right-edge speed
	of the classical contact process with infection
	rate $\lambda$.
\end{proposition}

The regenerative structure also yields a functional
central limit theorem for the fluctuations of the
right edge around its deterministic linear growth.

\begin{theorem}[Functional central limit theorem
	for the right edge]
	\label{teo_fclt}
	Fix $\lambda>\lambda_c(\infty)$ and let
	$\mu\geq\mu_0(\lambda)$. Let $v=v(\lambda,\mu)$
	be the asymptotic speed from
	Theorem~\ref{teo_lln}. Then there exists a
	deterministic constant
	\[
	\sigma^2=\sigma^2(\lambda,\mu)\in(0,\infty)
	\]
	such that, for every $\xi\in\mathcal X^{\mathrm{e}}$,
	under $\mathbb P^\xi_{\lambda,\mu}$,
	\[
	\left(
	\frac{r_{nt}-vnt}{\sqrt n}
	\right)_{t\geq0}
	\Longrightarrow
	\left(
	\sigma B_t
	\right)_{t\geq0},
	\qquad n\to\infty,
	\]
	in $D([0,\infty),\mathbb R)$ equipped with the
	Skorokhod $J_1$ topology, where $(B_t)_{t\geq0}$
	is a standard Brownian motion.
\end{theorem}

Beyond the motion of the infection front, we study
the configuration of the process as seen from its
right edge. The shielding construction also provides
a coupling mechanism: infinite infection paths form barriers that prevent differences in the initial configurations from influencing the evolution on the protected side.

For a nonempty configuration $\xi_t$ with a finite rightmost infected site, define
\[
\Xi_t(x):=\xi_t(r_t+x),
\qquad x\in\mathbb Z.
\]

\begin{theorem}[Invariant law seen from the right edge]
	\label{teo:edge-invariant-law}
	Fix $\lambda>\lambda_c(\infty)$ and let
	$\mu\geq\mu_0(\lambda)$. The CPTI seen from its
	right edge admits a unique invariant probability
	measure $\nu^{\mathrm{e}}$ on
	$\mathcal X^{\mathrm{e}}$.
	
	Moreover, for every $\xi\in\mathcal X^{\mathrm{e}}$,
	\begin{equation}
    \label{eq:weak-convergence-edge}
	\mathbb P^\xi_{\lambda,\mu}(\Xi_t\in\cdot)
	\Longrightarrow
	\nu^{\mathrm{e}},
	\qquad t\to\infty.
	\end{equation}
\end{theorem}

This result is the CPTI analogue of the edge convergence
theorem for the one-dimensional supercritical contact
process \cite{GalvesPresutti1987}. Its proof follows the
finite-window approach of \cite{CoxDurrettSchinazi1991},
using the shielding property to show that the distribution in every fixed window behind the edge becomes asymptotically independent of the initial configuration.

\begin{remark}\label{rem:uniform-edge-convergence}
The methods developed here can be extended to yield
a uniform version of the preceding convergence result.
More precisely, for every $0<c<v$ and every bounded
local function $f$,
\[
\lim_{t\to\infty}
\sup_{\substack{
A\subset\mathbb Z,\ 0<|A|<\infty\\
\operatorname{diam}(A)\le ct
}}
\left|
\mathbb E^A_{\lambda,\mu}
\left[f(\Xi_t)\mid\tau^A>t\right]
-\int f\,d\nu^{\mathrm e}
\right|=0.
\]
We refer to Remark~\ref{rem:uniform-edge-convergence-strategy}
below for a discussion of a possible strategy to
establish this uniform convergence.
\end{remark}

Let $\mathbf 0\in\{-1,0,1\}^{\mathbb Z}$ denote
the all-susceptible configuration,
\[
\mathbf 0(x)=0,
\qquad x\in\mathbb Z.
\]
This configuration is absorbing, and we write
$\delta_{\mathbf 0}$ for the corresponding
point mass.

The shielding mechanism also applies to coupling in
fixed spatial windows. By constructing barriers that
propagate in opposite directions, we can isolate
finite regions from differences in the surrounding
configurations and their infection histories. This
two-sided coupling provides the key ingredient for
our complete convergence theorem. 

\begin{theorem}[Complete convergence]
	\label{teo:complete-convergence}
	Fix $\lambda>\lambda_c(\infty)$ and let
	$\mu\geq\mu_0(\lambda)$. There exists a nontrivial
	translation-invariant invariant probability measure
	$\nu$ for the CPTI such that, for
	every nonempty $A\subset\mathbb Z$,
	\[
	\mathbb P_{\lambda,\mu}(\xi_t^A\in\cdot)
	\Longrightarrow
	\mathbb P_{\lambda,\mu}(\tau^A<\infty)\,
	\delta_{\mathbf 0}
	+
	\mathbb P_{\lambda,\mu}(\tau^A=\infty)\,
	\nu,
	\qquad t\to\infty.
	\]
	Moreover,
	\[
	\mathbb P_{\lambda,\mu}
	\bigl(
	\xi_t^A\in\cdot
	\mid
	\tau^A=\infty
	\bigr)
	\Longrightarrow
	\nu,
	\qquad t\to\infty.
	\]
\end{theorem}

\subsection{Graphical construction and auxiliary toy process}
\label{sec:toy}

The shielding construction relies on an auxiliary
process introduced in~\cite{AndjelSchinazi1996},
whose infection paths provide the space-time
barriers used throughout the paper. We first recall
the graphical construction of the CPTI and then
introduce this auxiliary process.

We construct the process from a Harris graphical representation on the whole time line. For each
ordered nearest-neighbor pair $(x,y)$, $|x-y|=1$, we put infection arrows from $x$ to $y$ at the times of a Poisson process of rate $\lambda$ on $\mathbb R$. For each site $x$, we put recovery marks at the times of a Poisson process of rate $1$ on $\mathbb R$. To each recovery mark $(x,t)$ we attach an independent exponential random variable $I(x,t)$ with parameter $\mu$. All Poisson processes and all variables $I(x,t)$ are independent.
 
If the initial configuration contains immune sites, then for each
$x$ with $\xi_0(x)=-1$ we assign an independent residual immunity time $I_x\sim \operatorname{Exp}(\mu)$,
independent of the Harris system and of the residual immunity times of the
other initially immune sites. At time $I_x$, the site $x$ becomes susceptible.

Given an initial configuration $\xi_0$, the process $(\xi_t)_{t\ge0}$ is obtained by reading the graphical construction after time $0$ in increasing time order.  At a recovery mark $(x,t)$, if $\xi_{t-}(x)=1$, the mark is \emph{effective}: we set $\xi_t(x)=-1$, and the site $x$ remains immune during the interval $[t,t+I(x,t))$.  At time $t+I(x,t)$, the site returns to state $0$.  If $\xi_{t-}(x)\neq1$, the recovery mark has no effect, and its attached interval is not activated.

We distinguish between two notions of infection path.
A~\emph{graphical arrow path} is a piecewise constant
nearest-neighbor space-time path that follows infection
arrows in their direction and does not cross recovery
marks at its current site.

Such a path is a \emph{legal CPTI infection path} if
it starts from an infected space-time point and never
enters an immune site. It may enter an already infected
site, in which case the paths coalesce. The infected
sites at any given time are precisely those reached
by legal infection paths from the initial infected set.
Unlike graphical arrow paths, legality depends on the
immunity intervals activated by the CPTI evolution.

We now introduce a stationary version of the auxiliary
process of~\cite{AndjelSchinazi1996}, which we call the
\emph{toy process}. The key idea is to make the blocking
mechanism independent of the infection history, thereby
recovering an additive structure.

The toy process uses the same infection arrows, recovery
marks, and attached exponential variables as the CPTI,
but treats immunity intervals differently. In the CPTI,
an interval is activated only when its recovery mark
corresponds to an actual recovery of an infected site.
In the toy process, every recovery mark creates a
blocking interval. Thus, a recovery mark at $(x,s)$
blocks the site $x$ during $[s,s+I(x,s))$, regardless
of whether the site is infected at time $s$.

Since the recovery marks are defined on all of
$\mathbb R$, blocking intervals born before time $0$
may still be active at that time. We write
\[
\beta_t(x)
:=
\#\{s\le t:\text{ there is a recovery mark at }(x,s)
\text{ and }s+I(x,s)>t\},
\qquad t\in\mathbb R,
\]
for the number of active blocking intervals at site
$x$ at time $t$. The site is toy-blocked precisely
when $\beta_t(x)>0$.

For each $x \in \mathbb Z$, the process $(\beta_t(x))_{t\in\mathbb R}$ is time-stationary. Moreover, for each $t \in \mathbb R$, the random variables $(\beta_t(x))_{x \in \mathbb Z}$ are i.i.d. with
\[
\beta_t(x)\sim
\operatorname{Poisson}(1/\mu)
\]
for each fixed $t$. Conditional on $\beta_0(x)=k$,
the $k$ blocking intervals active at time $0$ have
independent residual lifetimes with distribution
$\operatorname{Exp}(\mu)$. Equivalently, the blocking
field consists of independent immigration--death
processes, with immigration rate $1$ and death
rate $\mu$ per blocker.

Throughout the paper, we reserve the letter $\eta$
for the toy process. Given an initial infected set
$A\subseteq\mathbb Z$, let $\eta_t^A$ denote the
toy infected set at time $t$, started from $A$ at
time $0$. Sites in $A$ are initially infected
regardless of the value of $\beta_0$; active blockers
prevent only subsequent infections into blocked sites.

A toy infection is killed at a recovery mark at its
current site. An infection arrow from $x$ to $y$ at
time $t$ transmits infection if $x$ is infected
immediately before $t$ and $\beta_{t-}(y)=0$.
Equivalently, a \emph{toy infection path} is a graphical arrow path that never enters a blocked site after its starting time.
The path representation immediately implies that
the toy process is \emph{additive}.

The joint process $(\eta_t^A,\beta_t)$ is Markovian,
although $\eta_t^A$ alone is not, since its evolution
depends on the blocking field.
For survival probabilities, we condition the initial
infected seed to be unblocked; in particular, a
singleton started at $(0,0)$ is considered under
$\beta_0(0)=0$.

The crucial comparison with the CPTI is that every toy
infection path starting from a CPTI-infected site is
also a legal CPTI infection path, provided the two
processes use the same graphical construction and
all immunity intervals inherited from the initial
configuration are included among the toy-blocking
intervals. Thus, despite the lack of attractiveness of the CPTI,
the toy process provides infection paths that remain
valid under the original dynamics.

\begin{figure}[t]
    \centering
    \scalebox{0.9}{\input{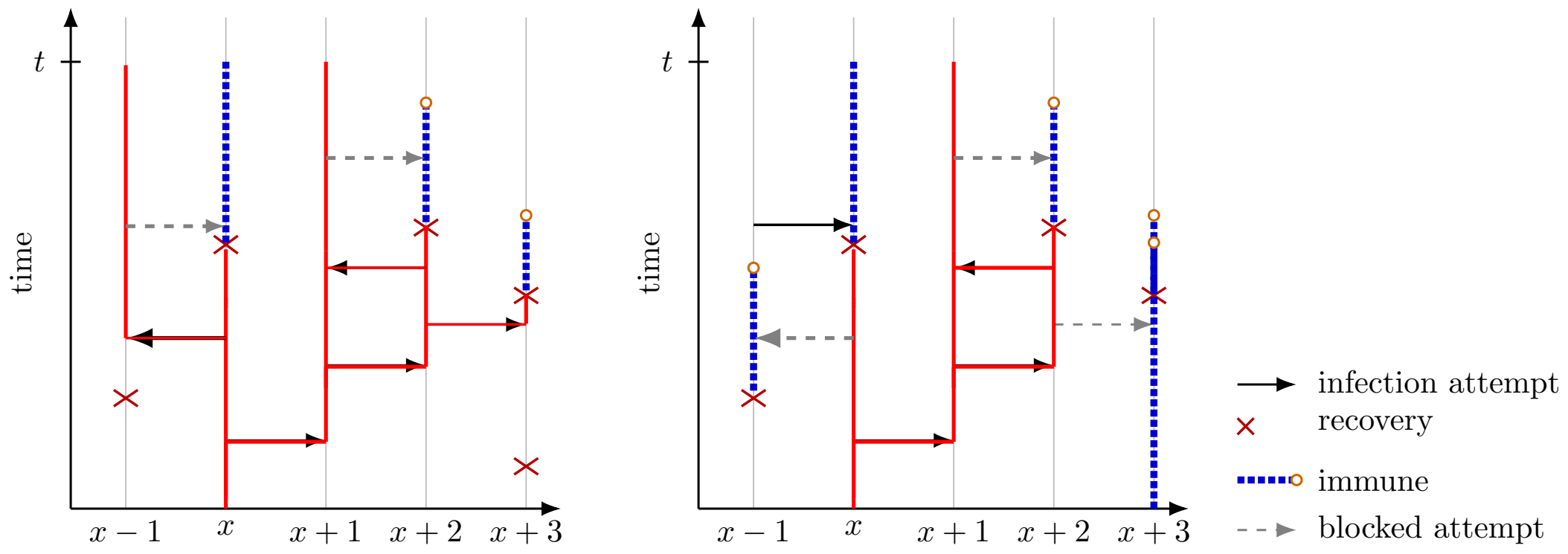}}
   \caption{Graphical construction of (a) the CPTI and
(b) the toy process. Solid arrows represent infection
attempts, red crosses recovery marks, dotted blue
segments immunity or blocking intervals, and dashed
arrows blocked attempts. Highlighted trajectories
are legal infection paths. The intervals at $x+3$
illustrate blockers inherited from the past and
overlapping blocking intervals.}
    \label{fig:graphical-construction} 
\end{figure}

Figure~\ref{fig:graphical-construction} illustrates
the CPTI and the toy process under a common
graphical construction.

\subsection{Overview of the approach}
\label{subsec:outline}

The autonomy of the toy-blocking field is crucial
for our approach. Unlike CPTI immunity, its evolution
does not depend on the infection history, yielding
an additive graphical representation. For sufficiently
large immunity-loss rates, recent blocking intervals
form a small perturbation of the classical contact
process, while older intervals constitute a sparse
field of obstacles. A block construction then yields survival and growth estimates for toy infection paths. These paths play a dual role: they transfer survival to the CPTI and provide the shielding barriers and quantitative estimates needed for regeneration.

The next difficulty is to turn survival into a
regenerative structure. The advance of the CPTI
right edge provides opportunities to explore fresh
graphical information and search for infinite toy
infection paths. By arranging these trials suitably,
we eventually obtain such a path, which can be used
to separate the evolution of the front from its
past. The selection of the leftmost infinite path is essential: its extremal character yields the one-sided shielding property needed to separate the relevant graphical information.

To understand this shielding mechanism, observe that
any influence from the left of the barrier must cross
it before reaching its right side. Since the barrier
is a legal infection path for all CPTI realizations
under consideration, its sites are infected in every
such realization. Consequently, differences in the
configuration or graphical history on the left cannot
propagate across the barrier to affect the evolution
on its right. Although this does not provide a global
comparison principle, it gives a one-sided protection
property that is sufficient for both regeneration
and coupling.

This separation allows us to replace the graphical history on the left of the barrier by a suitable completion without changing the subsequent evolution of the right edge. The resulting construction yields regeneration times with independent and identically distributed increments after the first cycle. The block estimates
yield the moment bounds needed to control both the
regeneration increments and the fluctuations within
each cycle, leading to the law of large numbers and
the functional central limit theorem.

The shielding mechanism also replaces monotone coupling in the study of long-time behavior. Common shielding paths force CPTI evolutions started from different configurations to agree in expanding spatial regions, despite the absence of attractiveness. This leads to
convergence to a unique invariant law for the
process seen from the edge and, by applying shielding
from both spatial directions, to complete convergence.

Finally, our extinction criterion follows from a
separate genealogical argument specific to one
dimension. The nearest-neighbor geometry allows
return infections to be assigned to previously
traversed infection edges, whose contributions
can be controlled by a branching-process comparison.
This yields a subcritical bound on the expected
number of infections in successive generations.

\subsection{Organization of the paper}

Section~\ref{sec:block-renormalization-survival}
establishes survival and growth properties for the
toy process, yielding Theorem~\ref{teo_survive}.
In Section~\ref{sec:regeneration-renewal}, we develop
the shielding and regeneration constructions and
establish the renewal structure.
Section~\ref{sec:lln-clt} derives the regenerative
limit theorems and the convergence of the edge speed,
using exponential moment bounds whose proof is
deferred to the Appendix~\ref{app:moment-estimates}.
Section~\ref{sec:shielding-asymptotic-laws}
establishes convergence to the invariant law and complete convergence
through shielding couplings.
Finally, we prove the extinction criterion in
Theorem~\ref{teo_dies} in Section~\ref{sec:extinction}.

\section{Survival  and growth of the toy process}
\label{sec:block-renormalization-survival}

We use a block construction to establish survival and
growth of the toy process and obtain the estimates needed
for regeneration. The key observation is that, for large
immunity-loss rates, recent blocking intervals produce
only a small perturbation of the classical contact process,
allowing us to approximate them by classical blocks. Although
blocking intervals originating in the distant past introduce
long-range dependence, they are sufficiently rare to be
controlled as sparse obstacles through a Peierls argument.

\subsection{The cutoff block construction}
\label{subsec:cutoff-block-construction}

The toy-blocking field of Section~\ref{sec:toy}
is stationary in time, so blocking intervals born
before time $0$ may affect the process afterward.
We first separate recent blocking intervals from
older ones.

Fix integers $H,\ell,K$, with $1\le \ell\le H$ and $K>H$, and fix a time scale $S>0$. For $j\in\mathbb Z$, set
\[
        I_j:=jK+[0,H].
\]
We say that a configuration \emph{contains a seed in $I_j$} if every subinterval of $I_j$ of length $\ell$ contains at least one infected site.  This condition is increasing in the infected set.

Let
\[
\mathcal L:=\{(j,n)\in\mathbb Z\times\mathbb Z_{\ge0}:
j+n\text{ is even}\}
\]
be the oriented lattice, with edges from $(j,n)$ to
$(j\pm1,n+1)$. A block $(j,n)$ propagates a seed
in $I_j$ at time $nS$ to seeds in $I_{j\pm1}$
at time $(n+1)S$.

Choose a finite interval $W_0\subset\mathbb Z$, and write $W_j:=jK+W_0$.  The space-time box used by block $(j,n)$ is
\[
        B_{j,n}:=W_j\times[nS,(n+1)S].
\]

Let $G^{\rm cl}_{j,n}$ be the event that every seed in
$I_j$ is connected to both neighboring seed intervals by classical contact paths with parameter $\lambda$ inside $B_{j,n}$. Note that the field of block indicators is finite-range dependent.
We choose $W_0$ large enough so that $G^{\rm cl}_{j,n}$ is determined by the graphical construction inside $B_{j,n}$.

Since $\lambda>\lambda_c(\infty)$, the classical contact
process admits a supercritical block construction.
For every $\varepsilon>0$, the parameters can be chosen
so that the block events $G^{\rm cl}_{j,n}$ satisfy
\[
\mathbb P(G^{\rm cl}_{j,n})\ge 1-\varepsilon.
\]

Fix $R>0$. Let $\mathcal D^{\le R}_{j,n}$ be the event
that an infection arrow in $B_{j,n}$ is suppressed by
a blocking interval of age at most $R$, born at a
recovery mark $(x,u)$ at the arrow's target $x$,
with $0\le t-u\le R$, where $t$ is the arrow time. We define $\mathcal D^{> R}_{j,n}$ analogously, but with blocking intervals with age at least $R$.

At a fixed site $x$ and time $t$, the expected number
of active blocking intervals of age at most $R$ is
\[
\int_0^R e^{-\mu a}\,da\le\frac1\mu.
\]
Since infection arrows enter each site at rate
$2\lambda$, we obtain
\[
\mathbf P_\mu(\mathcal D^{\le R}_{j,n})
\le\frac{2\lambda|W_0|S}{\mu}.
\]

Define the \emph{cutoff toy-good event} by
\[
\widetilde G^{{\rm toy},R}_{j,n}
:=
G^{\rm cl}_{j,n}\cap(\mathcal D^{\le R}_{j,n})^c.
\]
For fixed $R$, its indicator field is finite-range
dependent, and
\[
\mathbf P_\mu(\widetilde G^{{\rm toy},R}_{j,n})
\geq 1-\varepsilon-\frac{2\lambda|W_0|S}{\mu}.
\]
Thus its marginal probability can be made arbitrarily
close to one by choosing $\varepsilon$ small and then
$\mu$ large.

\subsection{Sparse old-tail obstacles}
\label{subsec:old-tail-obstacles}

We now control suppressions caused by blocking intervals
older than $R$. 
Rather than track suppressed arrows
directly, we introduce a larger obstacle event.
Let $\mathcal O_{j,n}$ be the event that some recovery
mark $(x,u)$, with $x\in W_j$ and $I(x,u)>R$, satisfies
\[
[u+R,u+I(x,u)]\cap[nS,(n+1)S]\neq\varnothing.
\]
Thus $\mathcal O_{j,n}$ records old blocking tails
intersecting the block. In particular,
$\mathcal D^{>R}_{j,n}\subseteq\mathcal O_{j,n}$,
so it suffices to control these obstacles.

We need the following uniform estimate.

\begin{lemma}[Sparse old-tail obstacles]
\label{lem:old-tail-obstacles}
Fix the block parameters $H,\ell,K,S,W_0$ and the
cutoff $R>0$. There exists $\rho_E=\rho_E(\mu,R)$,
with $\rho_E(\mu,R)\to0$ as $\mu\to\infty$ for fixed
$R$, such that, for every finite $A\subset\mathcal L$,
\[
\mathbf P_\mu\bigl(
\mathcal O_z\text{ occurs for every }z\in A
\bigr)
\leq \rho_E^{|A|}.
\]
\end{lemma}

\begin{proof}
We dominate the old tails by a Poisson stick process.
For each recovery mark $(x,u)$ with $I(x,u)>R$, set
\[
V:=u+R,\qquad L:=I(x,u)-R.
\]
By Poisson thinning and the memoryless property,
the points $(x,V)$ form a Poisson process of
intensity $e^{-\mu R}$ on each vertical line,
with independent lengths $L\sim\mathrm{Exp}(\mu)$.

A stick $(x,V,L)$ affects the set of block sites
\[
F(x,V,L)
:=
\left\{
(j,n)\in\mathcal L:
x\in W_j,\ 
[V,V+L]\cap[nS,(n+1)S]\neq\varnothing
\right\}.
\]
We call $F(x,V,L)$ the \emph{block footprint} of the stick.

A stick affects only a bounded number of block columns,
depending on $W_0$ and $K$, and at most $2+L/S$ time
levels, up to a change of constants. Since
$L\sim\mathrm{Exp}(\mu)$, its block footprint has an
exponential size tail, uniformly for large $\mu$.

Let $\nu_\mu$ be the intensity measure of the block
footprints. The stick intensity is $e^{-\mu R}$, and
footprint sizes have uniformly exponential tails.
Hence we can choose $\rho_E=\rho_E(\mu,R)\in(0,1)$,
with $\rho_E\to0$ as $\mu\to\infty$, such that
\[
\sup_{z\in\mathcal L}
\int_{\{F:z\in F\}}
\rho_E^{-|F|}\,\nu_\mu(dF)
\leq 1.
\]

We prove the claim by induction on $|A|$, the case
$A=\varnothing$ being trivial. Fix $z\in A$. If every
block in $A$ is hit, some stick has a footprint $F$
containing $z$, while all blocks in $A\setminus F$
must also be hit. By the Mecke formula and the
induction hypothesis,
\[
\begin{aligned}
\mathbf P_\mu\bigl(
\mathcal O_y\text{ for all }y\in A
\bigr)
&\leq
\int_{\{F:z\in F\}}
\rho_E^{|A\setminus F|}\,\nu_\mu(dF)\\
&\leq \rho_E^{|A|}
\int_{\{F:z\in F\}}
\rho_E^{-|F|}\,\nu_\mu(dF)
\leq \rho_E^{|A|}.
\end{aligned}
\]
This completes the induction.
\end{proof}

\subsection{Survival and growth after old-tail deletion}
\label{subsec:peierls-toy-survival}

We now combine the high-density, finite-range dependent
cutoff field with the sparse old-tail obstacles. Deleting
blocks hit by old tails recovers the full toy block
construction. 

Define the toy-good block after old-tail deletion by
\[
        G^{\rm toy}_{j,n}
        :=
        \widetilde G^{{\rm toy},R}_{j,n}\cap\mathcal O_{j,n}^c .
\]
For $z=(j,n)\in\mathcal L$, write
\[
        \widetilde\omega^R_z
        :=
        \Ind{\widetilde G^{{\rm toy},R}_{j,n}},
        \qquad
        \omega^{\rm toy}_z
        :=
        \Ind{G^{\rm toy}_{j,n}} .
\]
We call $(\widetilde\omega^R_z)_{z\in\mathcal L}$ \emph{the cutoff toy-good block field}.  This is the field obtained after ignoring blocking intervals whose age is larger than $R$. For fixed $R$, the cutoff toy-good block field is finite-range dependent,
with a fixed dependence range, and its marginal probability tends to one
as $\varepsilon\downarrow0$ and then $\mu\to\infty$. Therefore, by~\cite{LiggettSchonmannStacey1997}, it can be coupled with an independent
oriented site percolation field $(\Pi_z)_{z\in\mathcal L}$ of parameter
$p$ such that
\[
\Pi_z\leq\widetilde\omega^R_z
\qquad\text{for every }z\in\mathcal L,
\]
where $p$ can be chosen arbitrarily close to one.

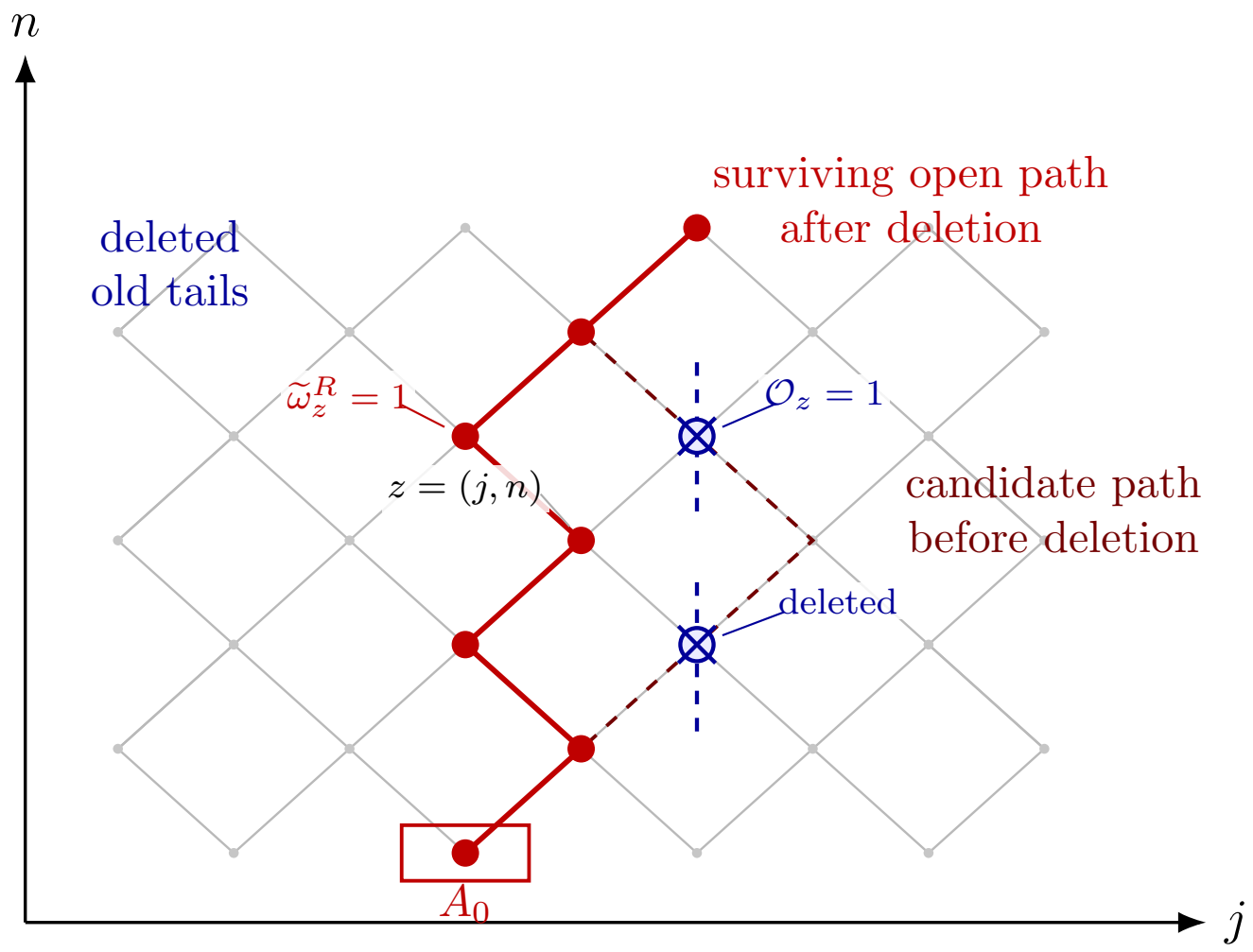
\begin{figure}[t]
    \centering
    \scalebox{0.90}{\centering
 \begin{tikzpicture}[x=0.9cm,y=0.9cm,>=Latex]

\tikzset{
    axis/.style={-{Latex[length=2.2mm]}, line width=0.7pt},
    lattice/.style={gray!55, line width=0.5pt},
    good/.style={
        circle,
        fill=red!75!black,
        draw=red!75!black,
        inner sep=2pt
    },
    openpath/.style={red!75!black, very thick},
    prepath/.style={red!45!black, thick, dashed},
    oldtail/.style={blue!60!black, dashed, line width=0.9pt},
    tailmark/.style={blue!60!black, line width=0.9pt},
    deletedblock/.style={
        circle,
        draw=blue!60!black,
        fill=blue!8,
        line width=0.9pt,
        inner sep=2.6pt
    },
    guide/.style={densely dashed, gray!70, line width=0.6pt},
    lab/.style={font=\small},
    blocklabel/.style={
        font=\scriptsize,
        fill=white,
        fill opacity=0.85,
        text opacity=1,
        inner sep=1.2pt
    },
    labelarrow/.style={
    gray!65,
    line width=0.45pt
}
}


\draw[axis] (-4.8,0) -- (-4.8,7.5) node[above] {$n$};
\draw[axis] (-4.8,0) -- (5.4,0) node[right] {$j$};


\foreach \x in {-3,-1,1,3}{
    \draw[lattice] (\x,0.6) -- ({\x-1},1.5);
    \draw[lattice] (\x,0.6) -- ({\x+1},1.5);
}

\foreach \x/\a/\b in {
    -4/-3/-3,
    -2/-3/-1,
     0/-1/1,
     2/1/3,
     4/3/3
}{
    \draw[lattice] (\x,1.5) -- (\a,2.4);
    \draw[lattice] (\x,1.5) -- (\b,2.4);
}

\foreach \x in {-3,-1,1,3}{
    \draw[lattice] (\x,2.4) -- ({\x-1},3.3);
    \draw[lattice] (\x,2.4) -- ({\x+1},3.3);
}

\foreach \x/\a/\b in {
    -4/-3/-3,
    -2/-3/-1,
     0/-1/1,
     2/1/3,
     4/3/3
}{
    \draw[lattice] (\x,3.3) -- (\a,4.2);
    \draw[lattice] (\x,3.3) -- (\b,4.2);
}

\foreach \x in {-3,-1,1,3}{
    \draw[lattice] (\x,4.2) -- ({\x-1},5.1);
    \draw[lattice] (\x,4.2) -- ({\x+1},5.1);
}

\foreach \x/\a/\b in {
    -4/-3/-3,
    -2/-3/-1,
     0/-1/1,
     2/1/3,
     4/3/3
}{
    \draw[lattice] (\x,5.1) -- (\a,6.0);
    \draw[lattice] (\x,5.1) -- (\b,6.0);
}

\foreach \x/\y in {
    -3/0.6,-1/0.6,1/0.6,3/0.6,
    -4/1.5,-2/1.5,0/1.5,2/1.5,4/1.5,
    -3/2.4,-1/2.4,1/2.4,3/2.4,
    -4/3.3,-2/3.3,0/3.3,2/3.3,4/3.3,
    -3/4.2,-1/4.2,1/4.2,3/4.2,
    -4/5.1,-2/5.1,0/5.1,2/5.1,4/5.1,
    -3/6.0,-1/6.0,1/6.0,3/6.0
}{
    \fill[gray!45] (\x,\y) circle (1.1pt);
}


\draw[red!75!black, thick] (-1.55,0.36) rectangle (-0.45,0.84);
\node[red!75!black, lab, below] at (-1,0.47) {$A_0$};


\draw[prepath]
    (-1,0.6) -- (0,1.5) -- (1,2.4) -- (2,3.3)
    -- (1,4.2) -- (0,5.1) -- (1,6.0);

\node[red!45!black, lab, align=center] at (4.08,3.57)
{candidate path\\before deletion};


\draw[oldtail] (1,1.65) -- (1,3.05);

\draw[oldtail] (1,3.55) -- (1,4.85);

\node[deletedblock] at (1,2.4) {};
\draw[tailmark] (0.84,2.24) -- (1.16,2.56);
\draw[tailmark] (0.84,2.56) -- (1.16,2.24);

\node[deletedblock] at (1,4.2) {};
\draw[tailmark] (0.84,4.04) -- (1.16,4.36);
\draw[tailmark] (0.84,4.36) -- (1.16,4.04);

\node[blue!60!black, lab, align=center] at (-3.55,5.70)
{deleted\\old tails};


\foreach \x/\y in {
    -1/0.6,0/1.5,-1/2.4,0/3.3,-1/4.2,0/5.1,1/6.0
}{
    \node[good] at (\x,\y) {};
}

\draw[openpath]
    (-1,0.6) -- (0,1.5) -- (-1,2.4) -- (0,3.3)
    -- (-1,4.2) -- (0,5.1) -- (1,6.0);

\node[red!75!black, lab, align=center] at (2.85,6.25)
{surviving open path\\after deletion};


\node[blocklabel] at (-1,3.75) {$z=(j,n)$};
\draw[labelarrow] (-0.35,3.70) -- (-0.10,3.42);

\node[red!75!black, blocklabel] at (-2.01,4.54)
{$\widetilde\omega^R_z=1$};
\draw[labelarrow, red!75!black] (-1.55,4.47) -- (-1.18,4.28);

\node[blue!60!black, blocklabel] at (2.22,2.78)
{deleted};
\draw[labelarrow, blue!60!black] (1.75,2.68) -- (1.22,2.48);

\node[blue!60!black, blocklabel] at (2.10,4.56)
{$\mathcal O_z=1$};
\draw[labelarrow, blue!60!black] (1.68,4.48) -- (1.22,4.28);

\end{tikzpicture}}
    \caption{Block renormalization with old-tail deletion.
The dashed red path is a candidate before deletion.
The blue segments represent old-tail obstacles
$\mathcal O_z$, which delete the crossed blocks.
The solid red path survives deletion.}
    \label{fig:old-tails-deletion}
\end{figure}

The next proposition says that deleting the sparse old-tail obstacles does not destroy the supercritical Peierls estimates, see Figure~\ref{fig:old-tails-deletion}.

\begin{proposition}
\label{prop:peierls-after-old-tail-deletion}
Suppose the cutoff field $\widetilde{\omega}^{R}$
dominates independent oriented site percolation with
parameter $p$, and set $q:=1-p$.
If $p$ is sufficiently close to one and $\rho_E$
sufficiently small, then the following estimates hold.

For every finite nonempty $A_0\subset\mathcal L$, let
\[
S^{\mathrm{blk}}(A_0)
:=\{\text{there is an infinite oriented path of
$G^{\mathrm{toy}}$-blocks from $A_0$}\}.
\]
Then
\[
\mathbb P_\mu(S^{\mathrm{blk}}(A_0))>0.
\]

Let $\mathcal C_n(A_0)$ be the set of block positions
reached at height $n$ by oriented paths of
$G^{\mathrm{toy}}$-blocks from $A_0$, and define
\[
H^{\mathrm{blk}}(A_0)
:=\sup\{n:\mathcal C_n(A_0)\neq\varnothing\}.
\]
Then there exist $C,c>0$ such that, for all $n\ge1$,
\[
\mathbb P_\mu\bigl(
H^{\mathrm{blk}}(A_0)\ge n,\,
(S^{\mathrm{blk}}(A_0))^c
\bigr)\le Ce^{-cn}.
\]

Define
\[
R_n^{\mathrm{blk}}(A_0)
:=\sup\mathcal C_n(A_0),
\]
with $\sup\varnothing=-\infty$.
Then, for some $\kappa,C,c>0$ and all $n\ge1$,
\[
\mathbb P_\mu\bigl(
S^{\mathrm{blk}}(A_0),\,
R_n^{\mathrm{blk}}(A_0)\le\max A_0+\kappa n
\bigr)\le Ce^{-cn}.
\]

On $S^{\mathrm{blk}}(A_0)$, let
$\Gamma^{\mathrm{blk}}(A_0)
=(\Gamma_n^{\mathrm{blk}}(A_0),n)_{n\ge0}$
be the leftmost infinite oriented path from $A_0$.
Possibly decreasing $\kappa$, there exist $C,c>0$
such that, for all $n\ge1$,
\[
\mathbb P_\mu\bigl(
S^{\mathrm{blk}}(A_0),\,
\Gamma_n^{\mathrm{blk}}(A_0)
\ge\min A_0-\kappa n
\bigr)\le Ce^{-cn}.
\]
\end{proposition}

\begin{proof}
Call a block $z$ \emph{bad} if $\Pi_z=0$ or
$\mathcal O_z$ occurs. Since
$\Pi_z\le\widetilde\omega^R_z$, every non-bad
block is a $G^{\rm toy}$-block. It therefore
suffices to consider oriented paths of non-bad blocks.

Let $A\subset\mathcal L$ be finite. If every
block in $A$ is bad, then for some $C\subseteq A$,
we have that $\mathcal{O}_z$ occurs for all $z \in C$ and $\Pi_z=0$
for every $z\in A\setminus C$. By
Lemma~\ref{lem:old-tail-obstacles} and the
product law of $(\Pi_z)$,
\[
\begin{aligned}
\mathbf P_\mu(A\text{ is entirely bad})
&\leq
\sum_{C\subseteq A}
\min\{\rho_E^{|C|},q^{|A|-|C|}\}\\
&\leq
(\sqrt{\rho_E}+\sqrt q)^{|A|}.
\end{aligned}
\]
Writing $\delta:=\sqrt{\rho_E}+\sqrt q$, we obtain
\begin{equation}
\mathbf P_\mu(A\text{ is entirely bad})
\leq\delta^{|A|}.
\label{eq:bad-block-bound}
\end{equation}
By choosing $p$ close to one and then $\mu$
large, $\delta$ can be made arbitrarily small.

The remaining argument uses standard planar
Peierls estimates for oriented percolation;
see~\cite{Durrett1984}. The essential point is
that~\eqref{eq:bad-block-bound} applies to every
prescribed dual barrier, without requiring
independence of the bad-block indicators.

Absence of an infinite oriented path from $A_0$
forces a finite connected dual barrier of bad
blocks separating $A_0$ from infinity. The number
of such barriers of size $m$ is at most $C_0^m$.
Hence, for $\delta$ sufficiently small,
\[
\mathbf P_\mu\bigl((S^{\rm blk}(A_0))^c\bigr)
\leq\sum_{m\geq m_0}C_0^m\delta^m<1,
\]
where $m_0\geq1$ is the minimal barrier size.
Thus $\mathbf P_\mu(S^{\rm blk}(A_0))>0$.

It remains to prove the last assertion. If the reached cluster has height at least $n$ but is finite, its enclosing bad barrier has
size at least $c_0n$. By~\eqref{eq:bad-block-bound},
\[
\begin{aligned}
&\mathbf P_\mu\bigl(
H^{\rm blk}(A_0)\geq n,\,
(S^{\rm blk}(A_0))^c
\bigr)\\
&\qquad\leq
\sum_{m\geq c_0n}C_0^m\delta^m
\leq Ce^{-cn}.
\end{aligned}
\]

The same dual-barrier argument in a rightward
cone gives the lower-speed estimate. For
$\kappa>0$ sufficiently small, on
$S^{\rm blk}(A_0)$ the event
$R_n^{\rm blk}(A_0)\leq\max A_0+\kappa n$
forces a connected bad barrier of size at
least $c_1n$, separating the reached cluster
from the right side of the cone. Hence, by
\eqref{eq:bad-block-bound},
\begin{equation}
\mathbf P_\mu\bigl(
S^{\rm blk}(A_0),\,
R_n^{\rm blk}(A_0)\leq\max A_0+\kappa n
\bigr)\leq Ce^{-cn}.
\label{eq:peierls-right-cone}
\end{equation}

By reflection, the same argument applies to
the leftmost infinite path. On
$S^{\rm blk}(A_0)$, the event
$\Gamma_n^{\rm blk}(A_0)\geq\min A_0-\kappa n$
forces a connected bad barrier of size at
least $c_2n$ along the leftward cone. Thus
\begin{equation}
\mathbf P_\mu\bigl(
S^{\rm blk}(A_0),\,
\Gamma_n^{\rm blk}(A_0)
\geq\min A_0-\kappa n
\bigr)
\leq Ce^{-cn}.
\label{eq:peierls-left-cone}
\end{equation}
\end{proof}

By additivity and translation invariance, positive toy survival
from $[0,H]$ implies positive survival from $\{0\}$.
For the regeneration construction, define
\[
\theta^{\rm toy}:=
\mathbf P_\mu\left(
\text{the toy process started from $\{0\}$ survives}
\,\middle|\,\beta_0(0)=0
\right)>0.
\]
Indeed, conditioning on $\beta_0(0)=0$ removes only
blocking intervals initially active at the origin
and cannot destroy any toy infection path.

\begin{corollary}[One-sided toy survival]
\label{cor:one-sided-toy-survival}
For $\mu$ sufficiently large,
\[
\theta^{\rm toy}_+
:=
\mathbf P_\mu\left(
\begin{array}{c}
\text{there exists an infinite toy infection path from $(0,0)$}\\
\text{contained in $[0,\infty)\times[0,\infty)$}
\end{array}
\,\middle|\,
\beta_0(0)=0
\right)>0.
\]
\end{corollary}

\begin{proof}
Run the block construction on the right half of the
oriented lattice, shifting the initial block so that
all boxes lie in the nonnegative half-line. The cutoff
estimate, Lemma~\ref{lem:old-tail-obstacles}, and the
Peierls argument of
Proposition~\ref{prop:peierls-after-old-tail-deletion}
still apply, giving positive probability of an
infinite toy infection path from a finite seed.

A finite graphical event in the nonnegative half-line
produces such a seed from $(0,0)$. Conditional on
$\beta_0(0)=0$, this event requires only finitely many
arrows, the absence of finitely many recovery marks,
and the absence of blockers on finitely many
space-time segments, so it has positive probability.
Concatenating the two constructions proves the claim.
\end{proof}

\begin{proof}[Proof of Theorem~\ref{teo_survive}]
Fix a finite nonempty $A\subset\mathbb Z$ and $x\in A$.
Every toy infection path from $\{x\}$ is a legal CPTI
infection path from $A$. Hence
\[
\mathbb P_{\lambda,\mu}(\tau^A=\infty)
\geq
\mathbf P_\mu\bigl(
\text{the toy process from $\{x\}$ survives}
\bigr)>0.
\]
\end{proof}

\begin{remark}
\label{rem:classical-speed-blocks}
Fix $\varepsilon>0$ and $p_0<1$. By the classical
contact-process growth theorem, the block dimensions
can be chosen so that
\[
\mathbf P(G_{j,n}^{\rm cl})>p_0
\]
and one rightward block step has microscopic slope
greater than $v_{\rm CP}(\lambda)-\varepsilon$.

Moreover, the right-edge large-deviation estimate
for supercritical oriented percolation
\cite[Theorem~11.1]{Durrett1984} gives block speed
arbitrarily close to one when the percolation parameter
is sufficiently large. Thus an infinite path of
$G^{\rm cl}$-blocks following this speed produces
a microscopic infection path with asymptotic speed
at least $v_{\rm CP}(\lambda)-\varepsilon$.
\end{remark}

\begin{remark}
\label{rem:toy-unbounded-growth}
On toy survival, the process reaches arbitrarily large
levels almost surely. Indeed, by
Proposition~\ref{prop:peierls-after-old-tail-deletion},
failed block attempts have finite height almost surely,
while each new attempt has a uniformly positive probability
of producing an infinite block path. A standard restart
argument completes the proof.
\end{remark}

\section{Regeneration and renewal}
\label{sec:regeneration-renewal}

The aim of this section is to construct a regenerative
structure for the right edge. The key idea is to use
successive record times to search for a fresh graphical
environment, free of inherited immunities. Within such
an environment, we select an infinite toy infection path
that acts as a space-time barrier, shielding the new
evolution from the information accumulated in the past.

\subsection{Prepared trials and the one-sided law}
\label{subsec:prepared-walls-Q}

Successive record levels of the right edge provide
opportunities to explore previously unvisited regions.
At such levels, we seek a local configuration free of
immunity from which a toy infection path can be launched.
If this path survives, it provides a candidate barrier
for regeneration.

Fix an initial condition $\xi_0 \in \mathcal{X}^{\mathrm{e}}$.
As $\xi_0$ is fixed, we will omit the dependence on the
initial condition. Let $T_x:=\inf\{t\ge0:r_t=x\}$.
Fix $s_*>0$. On $\{T_x<\infty\}$, let
$\mathcal W_x$ be the event that, during
$[T_x,T_x+s_*]$,

\begin{itemize}
\item an infection arrow from $x$ to $x+1$ occurs
before the first recovery mark at $x$ after $T_x$;
\item no recovery mark occurs at $x+1$ during
$[T_x,T_x+s_*]$;
\item no infection arrows occur from $x+1$ to $x+2$ during
$[T_x,T_x+s_*]$;
\item $\beta_{T_x+s_*}(x+1)=0$.
\end{itemize}

Only after this preparation do we launch the auxiliary toy process from the singleton $\{x+1\}$.  The toy-blocking field is not reinitialized at the launch time: the launched process uses the stationary blocking field defined
in Section~\ref{sec:toy}, including any blocking intervals born before the
launch time that are still active.

Throughout the regeneration construction, after rooting a launched toy process at time $0$, let
\[
        \tau^{\rm toy}:=\inf\{t\ge0:\eta_t^{\{0\}}=\varnothing\}
\]
be the extinction time of the toy process launched from the singleton, with
the convention $\tau^{\rm toy}=\infty$ on survival, and set
\[
        S^{\rm toy}:=\{\tau^{\rm toy}=\infty\}.
\]

We next fix the reference post-preparation law used throughout the section.
Starting from the stationary graphical construction of Section~\ref{sec:toy},
condition on the preparation event $\mathcal W_x$, shift the prepared site
and launch time to $(0,0)$, and retain the rooted graphical environment from
which the toy process is launched, including the active toy-blocking
intervals crossing time $0$ away from the prepared site. We denote this
biased post-preparation law by $\widehat{\mathbf P}$. By construction,
\[
        \widehat{\mathbf P}\bigl(\beta_0(0)=0\bigr)=1.
\]
Set
\[
        \widehat\theta:=\widehat{\mathbf P}(S^{\rm toy}).
\]
Lemma~\ref{lem:positive-survival-trial-law} below shows that
$\widehat\theta>0$ for $\mu$ sufficiently large.

\begin{lemma}[Positive survival under the trial law]
\label{lem:positive-survival-trial-law}
Assume that $\mu$ is in the regime of
Corollary~\ref{cor:one-sided-toy-survival}. Then
\[
        \widehat\theta
        =\widehat{\mathbf P}(S^{\rm toy})
        \ge \theta^{\rm toy}_+>0.
\]
\end{lemma}

\begin{proof}
Consider the post-preparation environment under $\widehat{\mathbf P}$ and
root the prepared site and launch time at $(0,0)$. By construction,
$\beta_0(0)=0$.

To prove survival it is enough to restrict the launched toy process to the right half-line $[0,\infty)$.

Before the preparation, the graphical information in this half-line has not been inspected. Apart from the condition $\beta_0(0)=0$, the preparation imposes only a fixed finite set of local requirements. Its infection-arrow requirements do not alter the toy-blocking field, while its recovery-mark requirements are absence requirements and are therefore favorable for toy infection paths.

Consequently, the post-preparation toy environment in the right half-line
may be coupled with the stationary toy environment conditioned on
$\beta_0(0)=0$ so that every one-sided toy infection path in the latter
environment is also available under $\widehat{\mathbf P}$. Therefore
Corollary~\ref{cor:one-sided-toy-survival} gives
\[
        \widehat\theta
        =\widehat{\mathbf P}(S^{\rm toy})
        \ge \theta^{\rm toy}_+>0.
\]
\end{proof}

Let $\mathcal E$ denote the rooted post-preparation graphical environment.
On $S^{\rm toy}$, the auxiliary toy construction contains at least one
infinite admissible infection path issued from the prepared site.  Among all
such paths, we choose the leftmost infinite path, with ties resolved according
to a fixed measurable deterministic rule, and denote it by
\[
        \Gamma=\Gamma(\mathcal E).
\]
This distinguished path will serve as the regeneration barrier.

We call the following procedure a \emph{trial based at level $x$}.  At the
record time $T_x$, we attempt the preparation event $\mathcal W_x$.

If $\mathcal W_x$ fails, the trial is rejected and is resolved at time
$T_x+s_*$.  If $\mathcal W_x$ occurs, then at time $T_x+s_*$ the site $x+1$ is an infected, unblocked right edge. From this site we launch the auxiliary toy process in the retained toy-blocking field.

If the launched toy process realizes $S^{\rm toy}$, the trial is accepted.
Among all infinite admissible toy infection paths issued from $x+1$, we choose
the leftmost infinite path according to the deterministic measurable rule
fixed above and denote it by $\Gamma_x$.

If the launched toy process does not realize $S^{\rm toy}$, the trial is
rejected.  In this case $\tau_x^{\rm toy}<\infty$, where
$\tau_x^{\rm toy}$ denotes the extinction time of the toy process launched
from $x+1$, measured from the launch time.  Thus the failure is resolved when
the launched toy process becomes empty. In summary, a trial has three possible outcomes:
\begin{itemize}
        \item{preparation failure},
        \item{preparation success but toy extinction},
        \item{preparation success and toy survival}.
\end{itemize}
Only the last outcome is accepted. The two last possibilities are illustrated
in Figure~\ref{fig:trials}.

\begin{figure}[t]
    \centering
    \centering
\resizebox{\textwidth}{!}{
\begin{tikzpicture}[x=0.78cm,y=0.78cm,>=Latex]

\tikzset{
    axis/.style={-{Latex[length=2.2mm]}, line width=0.7pt},
    guide/.style={densely dashed, gray!70, line width=0.6pt},
    edge/.style={
        black,
        thick,
        decorate,
        decoration={random steps, segment length=5pt, amplitude=0.55pt}
    },
    prepsite/.style={circle, fill=red!75!black, inner sep=1.8pt},
    toy/.style={red!75!black, very thick, line cap=round, line join=round},
    toybranch/.style={red!75!black, thick, line cap=round, line join=round},
    endpoint/.style={circle, fill=red!75!black, inner sep=1.2pt},
    lab/.style={font=\small},
    title/.style={font=\bfseries\large}
}

\def\Ty{3.9}       
\def\Wy{5.8}       
\def\Ry{9.35}      

\begin{scope}[shift={(0,0)}]

\node[title] at (4.6,10.6) {(a) Rejected trial};

\draw[axis] (0,0) -- (0,9.9) node[above] {time};
\draw[axis] (0,0) -- (9.0,0) node[right] {space};

\coordinate (record) at (3.2,\Ty);
\coordinate (prep)   at (3.9,\Wy);

\draw[edge]
    (0.7,0.35)
    .. controls (0.85,1.10) and (1.10,1.65) .. (1.25,2.15)
    .. controls (1.35,2.55) and (1.75,2.95) .. (2.00,3.25)
    .. controls (2.25,3.55) and (2.75,3.80) .. (3.20,\Ty);

\node[lab] at (1.15,3.15) {$r_t$};

\draw[guide] (3.2,0) -- (3.2,\Ty);
\draw[guide] (0,\Ty) -- (3.2,\Ty);
\node[below] at (3.2,0) {$x$};
\node[left] at (0,\Ty) {$T_x$};

\draw[guide] (0,\Wy) -- (4.35,\Wy);
\node[left] at (0,\Wy) {$T_x+s_*$};

\draw[<->, thick] (3.45,\Ty) -- (3.45,\Wy);
\node[lab, right] at (3.45,{(\Ty+\Wy)/2}) {$s_*$};

\node[prepsite] at (prep) {};
\node[below] at (3.9,0) {$x+1$};

\node[red!75!black, lab, align=center] at (8.10,8.05)
{finite toy\\exploration};

\draw[toy]
    (3.90,\Wy)
    -- (3.90,6.12)
    -- (4.35,6.12)
    -- (4.35,6.74)
    -- (4.80,6.74)
    -- (4.80,7.08)
    -- (5.25,7.08)
    -- (5.25,7.92);

\draw[toybranch]
    (4.35,6.42)
    -- (3.90,6.42)
    -- (3.90,7.26);

\draw[toybranch]
    (4.80,6.96)
    -- (4.35,6.96)
    -- (4.35,7.55);

\draw[toybranch]
    (5.25,7.36)
    -- (4.80,7.36)
    -- (4.80,7.84)
    -- (5.25,7.84)
    -- (5.25,8.28);

\draw[toybranch]
    (5.25,7.64)
    -- (5.70,7.64)
    -- (5.70,8.06)
    -- (6.15,8.06)
    -- (6.15,8.54)
    -- (6.60,8.54)
    -- (6.60,8.92);

\draw[toybranch]
    (6.15,8.36)
    -- (5.70,8.36)
    -- (5.70,8.72);

\draw[toybranch]
    (4.80,7.08)
    -- (4.80,7.42);

\node[endpoint] at (3.90,\Wy) {};
\node[endpoint] at (3.90,7.26) {};
\node[endpoint] at (4.35,7.55) {};
\node[endpoint] at (5.25,7.92) {};
\node[endpoint] at (5.25,8.28) {};
\node[endpoint] at (6.60,8.92) {};
\node[endpoint] at (5.70,8.72) {};
\node[endpoint] at (4.80,7.42) {};

\draw[blue!60!black, densely dotted, thick]
    (3.62,5.48)
    .. controls (3.33,6.22) and (3.40,7.42) .. (3.88,8.05)
    .. controls (4.22,8.52) and (4.85,8.72) .. (5.46,8.70)
    .. controls (5.95,9.18) and (6.64,9.35) .. (6.94,8.86)
    .. controls (7.24,8.28) and (6.82,7.58) .. (6.25,7.18)
    .. controls (6.02,6.46) and (5.72,5.82) .. (5.30,5.55)
    .. controls (4.70,5.40) and (4.08,5.38) .. cycle;

\draw[guide] (0,\Ry) -- (7.05,\Ry);
\node[left] at (0,\Ry) {$T_x+s_*+\tau_x^{\rm toy}$};

\node[lab, align=center] at (8.30,6.85)
{trial rejected;\\search continues};

\end{scope}

\begin{scope}[shift={(10.8,0)}]

\node[title] at (4.6,10.6) {(b) Accepted trial};

\draw[axis] (0,0) -- (0,9.9) node[above] {time};
\draw[axis] (0,0) -- (9.0,0) node[right] {space};

\coordinate (record) at (3.2,\Ty);
\coordinate (prep)   at (3.9,\Wy);

\draw[edge]
    (0.7,0.35)
    .. controls (0.85,1.10) and (1.10,1.65) .. (1.25,2.15)
    .. controls (1.35,2.55) and (1.75,2.95) .. (2.00,3.25)
    .. controls (2.25,3.55) and (2.75,3.80) .. (3.20,\Ty);

\node[lab] at (1.15,3.15) {$r_t$};

\draw[guide] (3.2,0) -- (3.2,\Ty);
\draw[guide] (0,\Ty) -- (3.2,\Ty);
\node[below] at (3.2,0) {$x$};
\node[left] at (0,\Ty) {$T_x$};

\draw[guide] (0,\Wy) -- (4.35,\Wy);
\node[left] at (-1.45,5.72) {$\tau_+=T_x+s_*$};

\draw[<->, thick] (3.45,\Ty) -- (3.45,\Wy);
\node[lab, right] at (3.45,{(\Ty+\Wy)/2}) {$s_*$};

\node[prepsite] at (prep) {};
\node[below] at (3.9,0) {$x+1$};

\draw[toy]
    (3.90,\Wy)
    -- (3.90,6.18)
    -- (4.35,6.18)
    -- (4.35,6.62)
    -- (4.80,6.62)
    -- (4.80,7.05)
    -- (4.35,7.05)
    -- (4.35,7.48)
    -- (4.80,7.48)
    -- (4.80,7.90)
    -- (5.25,7.90)
    -- (5.25,8.34)
    -- (4.80,8.34)
    -- (4.80,8.78)
    -- (5.25,8.78)
    -- (5.25,9.20)
    -- (5.70,9.20)
    -- (5.70,9.70);

\node[endpoint] at (3.90,\Wy) {};
\node[endpoint] at (4.35,6.18) {};
\node[endpoint] at (4.80,6.62) {};
\node[endpoint] at (4.35,7.05) {};
\node[endpoint] (node2) at (5.25,7.90) {};
\node[endpoint] (node1) at (4.80,8.34) {};
\node[endpoint] at (5.70,9.20) {};

\node[red!75!black, lab, right] at (5.78,8.72) {$\Gamma$};

\draw[toy,-{Latex[length=2.5mm]}]
    (5.70,9.25) -- (5.70,9.78);

\node[lab, align=center] at (1.95,8.25) {old information};
\node[lab, align=center] at (8.15,9.20) {future region};

\draw[edge]
    (3.90,\Wy)
    .. controls (4.70,6.25) and (5.45,6.80) .. (5.75,7.35)
    .. controls (6.05,7.70) and (6.30,8.20) .. (6.85,8.65);

\node[lab, align=center] at (8.05,6.40)
{future right edge};

\node[
    lab,
    align=center,
    draw=red!75!black,
    rounded corners,
    fill=red!4,
    inner sep=3pt
] at (7.55,4.05)
{$\Gamma$ is the\\regeneration barrier};

\draw[red!75!black,-{Latex[length=2mm]}]
    (7.32,4.70) -- (node2.south east);

\end{scope}

\end{tikzpicture}%
}
    \caption{Regeneration trials after successful preparation.
In (a), the toy process becomes extinct and the trial
is rejected. In (b), it survives, and its leftmost
infinite infection path $\Gamma$ becomes the
regeneration barrier.}
    \label{fig:trials}
\end{figure}
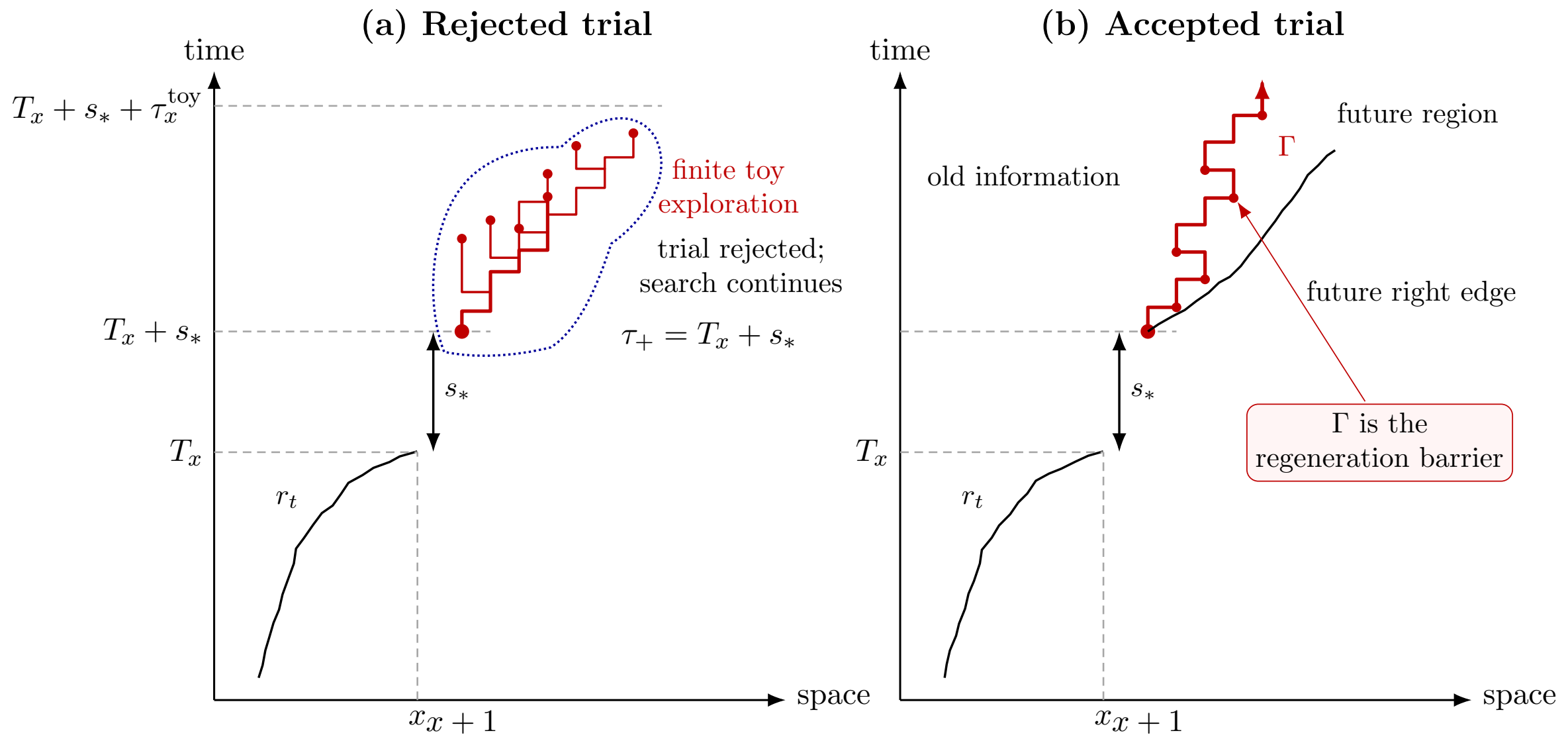 The finite graphical pattern entering $\mathcal W_x$ involves only a fixed
space-time box. The additional requirement $\beta_{T_x+s_*}(x+1)=0$ is imposed at the prepared site, whose toy-blocking history has not been
inspected before the trial. By stationarity of the toy-blocking field,
\[
        \mathbf P_\mu\bigl(\beta_0(0)=0\bigr)
        =e^{-1/\mu}>0.
\]
Since the remaining requirements in $\mathcal W_x$ form a fixed finite
graphical pattern compatible with this event, there exists
$q_{\rm prep}>0$ such that, uniformly over the information revealed before
the trial,
\begin{equation} \label{eq:conditional-lower-bound}
        \mathbf P_\mu(\mathcal W_x\mid\mathcal G_x)
        \ge q_{\rm prep}
        \qquad\text{on }\{T_x<\infty\},
\end{equation}
where $\mathcal G_x$ denotes the information revealed before the preparation
at level $x$ is performed.

By Lemmas~\ref{lem:positive-survival-trial-law} and~\ref{lem:freshness-adaptive-trials} below, once the preparation succeeds
the post-preparation environment has law $\widehat{\mathbf P}$ and the toy
survival probability is $\widehat\theta>0$. Hence every trial is accepted
with conditional probability at least
\[
        p_{\rm br}:=q_{\rm prep}\widehat\theta>0.
\]
An accepted barrier does not regenerate the full CPTI configuration. It
regenerates only the future right edge. We therefore introduce a one-sided
law $Q$. If $\mathcal E$ has law $\widehat{\mathbf P}$, condition on
$S^{\rm toy}$, select the leftmost infinite toy path $\Gamma$, and retain the
marked graphical construction on and to the right of $\Gamma$. Writing
$\Phi(\mathcal E)$ for the resulting rooted one-sided state, let
$ \mathsf\Omega_+$ denote the space of such rooted one-sided states and define
\[
        Q:=\mathcal L_{\widehat{\mathbf P}}
        \bigl(\Phi(\mathcal E)\mid S^{\rm toy}\bigr).
\]
We write $\mathbf P_Q$, $\mathbf E_Q$, and $\mathcal L_Q$ for probability,
expectation, and law when the initial rooted one-sided state has distribution
$Q$.

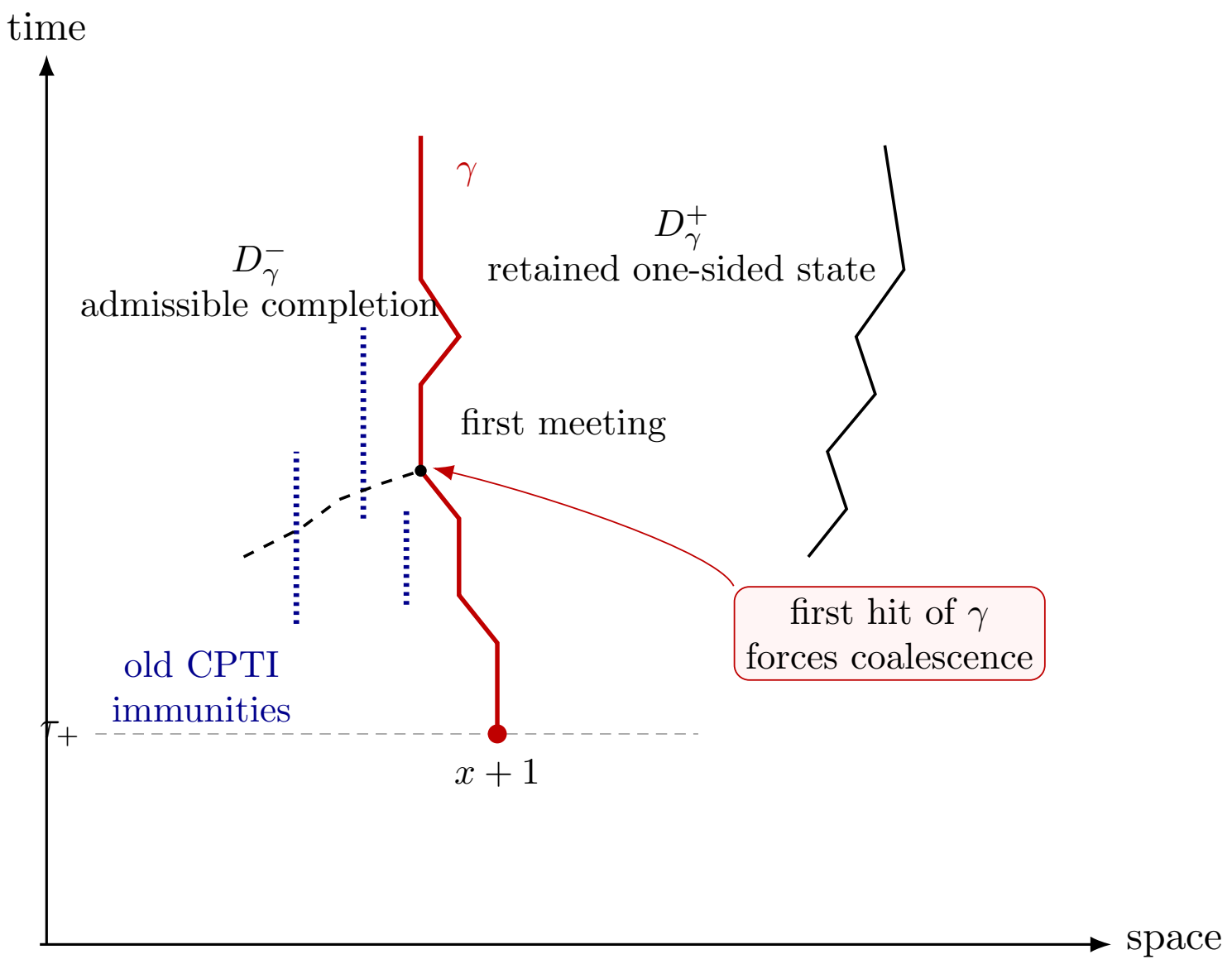
\begin{figure}[t]
\centering
\scalebox{0.85}{\centering
\begin{tikzpicture}[x=0.9cm,y=0.9cm,>=Latex]

\tikzset{
    axis/.style={-{Latex[length=2mm]}, line width=0.7pt},
    prepsite/.style={circle, fill=red!75!black, inner sep=1.8pt},
    barrier/.style={red!75!black, very thick},
    infpath/.style={black, thick},
    immune/.style={
        blue!55!black,
        line width=1.5pt,
        dash pattern=on 1pt off 1.4pt
    },
    guide/.style={densely dashed, gray!70},
    lab/.style={font=\small}
}

\draw[axis] (-0.51,0) -- (-0.51,9.3) node[above] {time};
\draw[axis] (-0.58,0) -- (10.62,0) node[right] {space};

\def\tauy{2.2}

\draw[guide] (0,\tauy) -- (6.3,\tauy);
\node[lab, left] at (0,\tauy) {$\tau_+$};

\node[prepsite] at (4.2,\tauy) {};
\node[lab, below=3pt] at (4.2,\tauy) {$x+1$};

\draw[immune] (2.1,3.35) -- (2.1,5.15);
\draw[immune] (2.8,4.45) -- (2.8,6.45);
\draw[immune] (3.25,3.55) -- (3.25,4.55);

\node[
    lab,
    blue!55!black,
    align=center
] at (1.11,2.7)
{old CPTI\\immunities};

\draw[barrier]
    (4.2,\tauy)
    -- (4.2,3.15)
    -- (3.8,3.65)
    -- (3.8,4.45)
    -- (3.4,4.95)
    -- (3.4,5.85)
    -- (3.8,6.35)
    -- (3.4,6.95)
    -- (3.4,8.45);

\node[
    lab,
    red!75!black,
    anchor=west
] at (3.62,8.05) {$\gamma$};

\node[
    lab,
    align=center
] at (1.72,6.92)
{$D_\gamma^-$\\admissible completion};

\node[
    lab,
    align=center
] at (6.13,7.33)
{$D_\gamma^+$\\retained one-sided state};

\draw[infpath, dashed]
    (1.55,4.05)
    -- (2.15,4.35)
    -- (2.55,4.65)
    -- (3.0,4.82)
    -- (3.4,4.95);

\fill (3.4,4.95) circle (1.6pt);

\node[
    lab,
    anchor=south east
] at (6.12,5.1)
{first meeting};

\node[
    lab,
    align=center,
    draw=red!75!black,
    rounded corners,
    fill=red!4,
    inner sep=3pt
] (coal) at (8.3,3.25)
{first hit of $\gamma$\\forces coalescence};

\draw[
    red!75!black,
    -{Latex[length=2mm]}
]
    (coal.north west)
    .. controls (6.5,4.1) and (4.7,4.7) ..
    (3.52,4.98);

\draw[infpath]
    (7.45,4.05)
    -- (7.85,4.55)
    -- (7.65,5.15)
    -- (8.15,5.75)
    -- (7.95,6.35)
    -- (8.45,7.05)
    -- (8.25,8.35);

\end{tikzpicture}}
\caption{Shielding mechanism. The infection path
$\gamma$ separates the retained one-sided graphical
state from an admissible completion on its left.
Any infection path entering the retained region
from the left must first meet $\gamma$, where
the paths coalesce.}
\label{fig:barrier-property}
\end{figure}

The future right edge can be read from a state with law $Q$ because the
distinguished leftmost toy path is a CPTI infection path and will act as a
barrier.
We next isolate the deterministic shielding property of an accepted barrier. The geometric mechanism is illustrated in
Figure~\ref{fig:barrier-property}.

For a càdlàg nearest-neighbor path $\gamma$, write $D_\gamma^+$ for the
closed space-time region on and to the right of $\gamma$, including the
horizontal edge traversed whenever $\gamma$ jumps, and write $D_\gamma^-$
for its strict-left complement.

We regard the marked graphical construction in $D_\gamma^+$, together with the distinguished path $\gamma$ and the CPTI state on and to the right of $\gamma(0)$ at time $0$, as a one-sided graphical state.

We now introduce the notion of \emph{completion} of a one-sided state. Informally, an admissible completion is any way of filling in the discarded
strict-left graphical history that is compatible with the retained one-sided
state and does not destroy the distinguished path $\gamma$ as a legal CPTI
infection path.

\begin{definition}An admissible completion of a one-sided state is any extension to
a full marked graphical construction on $\mathbb Z\times[0,\infty)$, together
with an initial CPTI configuration on the strict-left side, such that:
\begin{enumerate}
    \item the extension agrees with the retained graphical construction in
    $D_\gamma^+$ and with the retained CPTI state at time $0$ on and to the
    right of $\gamma(0)$; and
    \item in the resulting CPTI construction, $\gamma$ is a legal infection
    path starting from $\gamma(0)$ at time $0$.
\end{enumerate}
\end{definition}
Thus an admissible completion may be arbitrary on the discarded strict-left
side, subject only to compatibility with the retained one-sided state and to
the requirement that the distinguished path remain a legal CPTI infection
path.

\begin{lemma}[Shielding by an infection path]
\label{lem:toy-path-shielding}
Fix a one-sided graphical state consisting of an infinite path $\gamma$
together with the marked graphical construction in $D_\gamma^+$. Let $\xi^{(1)}$ and $\xi^{(2)}$ be the CPTI constructions obtained from two
admissible completions of the same one-sided graphical state.

Then, for every $t\ge0$,
\begin{equation}
\label{eq:shielding-configurations}
        \xi_t^{(1)}(x)=\xi_t^{(2)}(x)
        \qquad\text{for every }x\ge\gamma(t).
\end{equation}
Consequently their right edges agree for all $t\ge0$.
\end{lemma}

\begin{proof}
We first work in a finite space-time box and argue successively through its
graphical event times.

At time $0$ the two constructions agree on and to the right of
$\gamma(0)$.  Suppose that this agreement holds immediately before the next
relevant graphical event.

Any graphical event occurring strictly to the right of $\gamma$ acts
identically in the two constructions, since the marked graphical information
in $D_\gamma^+$ is common to both.  An infection arrow entering the current
site of $\gamma$ from its strict left cannot create a discrepancy, because
the site occupied by $\gamma$ is already infected in both constructions.

A recovery mark cannot occur on a vertical segment traversed by $\gamma$,
since $\gamma$ is a legal infection path in both constructions.

We now record the only model-specific input in the shielding argument.
The distinguished barrier $\gamma$ is a toy infection path, and under the
natural coupling every toy infection path is a legal CPTI infection path.
In particular, $\gamma$ cannot encounter an active CPTI immunity: if a site
visited by $\gamma$ were immune at the relevant time, $\gamma$ would not be
a legal CPTI infection path and hence could not be the barrier. Thus every
infection arrow traversed by $\gamma$ is legal for the CPTI.

Suppose that at time $u$ the path $\gamma$ jumps from $z$ to $z+1$.
Immediately before the jump, the two constructions agree at every site
$x\ge z$. Since the infection arrow traversed by $\gamma$ belongs to the
common retained graphical construction, it acts identically in both
constructions. Hence the two configurations still agree at every site
$x\ge z+1$, which is exactly the region on and to the right of $\gamma$
after the jump.

Suppose instead that at time $u$ the path $\gamma$ jumps from $z+1$ to $z$.
Immediately before the jump, the source $z+1$ is infected in both
constructions because it lies on $\gamma$.  Since both strict-left histories
are admissible completions, $\gamma$ is a legal infection path in each
construction.  In particular, the target site $z$ is available for the
arrow traversed by $\gamma$ in both constructions. Hence the arrow infects
$z$ in both constructions, or has no effect if $z$ is already infected.
Agreement therefore extends to the newly incorporated boundary site $z$.

Thus agreement on and to the right of the current position of $\gamma$ is
preserved after every graphical event.  Induction proves
\eqref{eq:shielding-configurations} in every finite space-time box.
Locality of the Harris construction, followed by an increasing-box
argument, gives the assertion for every finite time interval, and hence for
all $t\ge0$.

Finally, $\gamma(t)$ is infected in both constructions.  Since the two
configurations agree at every site on and to its right, their rightmost
infected sites coincide.
\end{proof}

\subsection{Construction of one cycle from \texorpdfstring{$Q$}{Q}}
\label{subsec:one-cycle-from-Q}

We now construct one regeneration cycle starting from the one-sided law $Q$. The rooted state is shifted so that the prepared site of the already accepted barrier is at the space-time point $(0,0)$.  Complete the graphical history on the strict-left side in any way compatible with the distinguished barrier.  In these coordinates the future CPTI right edge is denoted by $r_t$, with $r_0=0$; Lemma~\ref{lem:toy-path-shielding} shows that its future evolution is independent of the chosen compatible completion.  We write
\[
        \bar r_t:=\sup_{0\le s\le t}r_s .
\]
By Remark~\ref{rem:toy-unbounded-growth}, on $S^{\rm toy}$ the toy process reaches arbitrarily large levels almost surely. Since every toy infection path is a legal CPTI infection path under the natural coupling, it follows that $r_t$ reaches arbitrarily large levels almost surely under $Q$.

A rejected trial may reveal graphical information
beyond the current right edge, since the associated
toy process can explore far to the right before
becoming extinct. To retain fresh graphical
information for subsequent trials, we must therefore
keep track not only of the advance of the CPTI
right edge, but also of the region explored by
previous trials.

\medskip
\noindent\textbf{The revealed frontier.}
During rejected trials we reveal graphical information in order to decide
whether the preparation succeeded and, if it did, to follow the launched
toy process until extinction. A finite toy process may explore graphical
information far to the right of its starting level before it dies.  Set
$Y_0:=0$ and let $Y_i$ denote the rightmost spatial point whose graphical
information has been revealed after the first $i$ rejected trials have been resolved, including the fixed local information needed to decide the
preparation events and all information inspected while running the
rejected toy processes up to extinction.

\medskip
\noindent\textbf{Candidate record levels.}
Set $\rho_0:=0$. Suppose that the first $i-1$ trials have been rejected and resolved at times $\rho_1,\ldots,\rho_{i-1}$.  The next candidate level is chosen to be the first integer strictly to the right of both the current running record and all previously revealed trial information:
\begin{equation}
\label{eq:one-cycle-target-level}
        a_i:=
        1+\max\{\bar r_{\rho_{i-1}},Y_{i-1}\}.
\end{equation}
Let
\[
T_i:=\inf\{t\ge \rho_{i-1}:r_t=a_i\}=T_{a_i},
\]
where the second equality uses the notation of Subsection~\ref{subsec:prepared-walls-Q}. By Remark~\ref{rem:toy-unbounded-growth}, $T_i<\infty$ almost surely under $Q$.

At time $T_i$, the half-line strictly to the right of $a_i$ has not yet been reached by the CPTI right edge. Hence it contains no active CPTI immunities created by the past.  Moreover, by the definition of $a_i$ in~\eqref{eq:one-cycle-target-level}, this half-line contains no graphical information revealed while resolving the previous rejected trials. Thus the trial at level $a_i$ is performed in previously uninspected graphical information.

\medskip
\noindent\textbf{The $i$-th trial.}
At time $T_i$, we perform the trial based at level $a_i$, in the sense of Subsection~\ref{subsec:prepared-walls-Q}. If this trial is rejected, let $\rho_i$ be its resolution time.
For the $i$-th trial, write $\tau_i^{\rm toy}:=\tau_{a_i}^{\rm toy}$. Thus $\rho_i=T_i+s_*$ in the case of preparation failure, while
\[
    \rho_i=T_i+s_*+\tau_i^{\rm toy}
\]
in the case where the preparation succeeds but the launched toy process
becomes extinct.

For a rejected trial, define the trial overshoot by
\[        O_i:=
        \left(\max\{\bar r_{\rho_i},Y_i\}-a_i\right)_+ .
\]
Then, on rejection of the $i$-th trial,
\begin{equation}
\label{eq:one-cycle-target-recursion}
        a_{i+1}
        =
        1+\max\{\bar r_{\rho_i},Y_i\}
        \le
        a_i+O_i+1 .
\end{equation}

If the trial based at $a_i$ is accepted, the cycle stops.  Define
\begin{equation}
\label{eq:one-cycle-first-accepted-index}
        N:=\inf\{i\ge1:\text{the trial based at }a_i\text{ is accepted}\}.
\end{equation}
The accepted time, accepted position, and accepted barrier are
\[        \tau_+:=T_N+s_*,
        \qquad
        X_+:=a_N+1,
        \qquad
        \Gamma_+:=\Gamma_{a_N}.
\]

Let $\mathcal G_i$ denote the information revealed before the preparation
at level $a_i$ is performed.  By~\eqref{eq:one-cycle-target-level}, no
graphical information strictly to the right of $a_i$ has been inspected
while resolving the previous rejected trials. 

The preceding construction ensures that each candidate trial is performed
beyond all graphical information revealed by the earlier rejected trials.
The following lemma makes precise the resulting freshness property of the
post-preparation environment under this adaptive exploration.

\begin{lemma}[Invariant trial law under adaptive exploration]
\label{lem:freshness-adaptive-trials}
Let $\mathcal K_i$ be the $\sigma$-field generated by all the information
inspected before the post-preparation environment $\mathcal E_i$ of the $i$-th
trial is examined, including the information used to resolve the preceding
rejected trials and to perform the $i$-th preparation. Thus $\mathcal G_i$
contains the information before the $i$-th preparation is performed, whereas
$\mathcal K_i$ also contains the information used to decide that preparation.

On the event $\mathcal W_{a_i}$ that the $i$-th preparation succeeds,
\begin{equation}
\label{eq:fresh-environment-kernel}
\mathcal L_{Q}
\bigl(\mathcal E_i\mid\mathcal K_i\bigr)
=
\widehat{\mathbf P}
\qquad\mathbf{P}_{\mu}\text{ - a.s.}
\end{equation}
In particular, the conditional law does not depend on $i$ or on the realized
information in $\mathcal K_i$.

Moreover, for every bounded measurable function $G$ on the one-sided state 
space, on $\mathcal W_{a_i}$,
\begin{equation}
\label{eq:fresh-accepted-kernel}
\mathbf E_{Q}\!\left[
\Ind{S_i^{\rm toy}}
G\bigl(\Phi(\mathcal E_i)\bigr)
\,\middle|\,\mathcal K_i
\right]
=
\widehat\theta Q(G)
\qquad\mathbf{P}_{\mu}\text{ - a.s.},
\end{equation}
where
\[
S_i^{\rm toy}:=\{\mathcal E_i\in S^{\rm toy}\}.
\]
Consequently, on $\mathcal W_{a_i}$,
\[
\mathbf P_{Q}
\bigl(S_i^{\rm toy}\mid\mathcal K_i\bigr)
=
\widehat\theta
\qquad\mathbf{P}_{\mu}\text{ - a.s.},
\]
and, conditional on $\mathcal K_i$, $\mathcal W_{a_i}$, and
$S_i^{\rm toy}$, the accepted one-sided state $\Phi(\mathcal E_i)$ has law
$Q$.
\end{lemma}

\begin{proof}
The point is that the adaptive exploration does not inspect the graphical
information from which the next post-preparation environment is read.

First consider a one-sided state with law $Q$.  By definition, it is obtained
by conditioning a post-preparation environment with law
$\widehat{\mathbf P}$ on $S^{\rm toy}$, selecting the leftmost infinite toy
infection path, and retaining the marked graphical construction on and to the
right of that path. A left-first priority exploration identifies the leftmost infinite
toy path by following leftmost continuations and backtracking whenever
an alternative dies out. On $S^{\rm toy}$, this exploration reveals
only the selected path and failed alternatives on its left, leaving
the strict-right graphical field uninspected. By
planarity, a path lying strictly to the left of the selected infinite path
cannot cross it without first meeting it.  Thus neither toy survival nor the
choice of the leftmost surviving path requires inspection of the strict-right
marked graphical field.  Since the underlying marked Poisson construction is
a product field, this uninspected strict-right information retains its
original conditional law.  This is the unrevealed-field property built into
the one-sided law $Q$.

Now run the cycle exploration from $Q$.  Before the $i$-th trial, all
information revealed while resolving the preceding rejected trials lies at
or to the left of the revealed frontier $Y_{i-1}$.  By
\eqref{eq:one-cycle-target-level}, the candidate level $a_i$ is chosen
strictly to the right of both $Y_{i-1}$ and the running record
$\bar r_{\rho_{i-1}}$.  Thus, up to the first hitting time $T_i$, no
information strictly to the right of $a_i$ has been inspected by the cycle
exploration, and this half-line has not previously been reached by the CPTI
right edge. The preparation event then imposes exactly the conditioning used in the
definition of $\widehat{\mathbf P}$: the fixed local graphical pattern,
together with the requirement that the prepared launch site carry no active
toy blocker.

The toy-blocking field is part of the retained marked construction, including
blocking intervals born before the launch time. Nevertheless, on the
uninspected side its pre-launch marks have not been revealed. By the product
structure and space-time translation invariance of the stationary marked
construction, after translating the prepared site and launch time to $(0,0)$,
the conditional law of the resulting post-preparation environment is exactly
the fixed biased law $\widehat{\mathbf P}$. This proves
\eqref{eq:fresh-environment-kernel}.

It remains only to identify the accepted kernel.  By the definition of $Q$,
if $\mathcal E$ has law $\widehat{\mathbf P}$, then conditional on
$S^{\rm toy}$ the one-sided state $\Phi(\mathcal E)$ has law $Q$.
Therefore, for every bounded measurable $G$,
\[
\widehat{\mathbf P}\!\left[
\Ind{S^{\rm toy}}G\bigl(\Phi(\mathcal E)\bigr)
\right]
=
\widehat\theta Q(G).
\]
Combining this identity with \eqref{eq:fresh-environment-kernel} gives
\eqref{eq:fresh-accepted-kernel} on $\mathcal W_{a_i}$.  Taking $G\equiv1$
gives the stated conditional survival probability, and division by
$\widehat\theta$ gives the conditional law of $\Phi(\mathcal E_i)$ on
$\mathcal W_{a_i}\cap S_i^{\rm toy}$.
\end{proof}

 Conditional on a successful preparation, Lemma~\ref{lem:freshness-adaptive-trials} implies that the
post-preparation environment has law $\widehat{\mathbf P}$. Combining this
with the uniform preparation bound~\eqref{eq:conditional-lower-bound} gives
\[        \mathbf P_{Q}\left(
        \text{the trial based at }a_i\text{ is accepted}
        \mid \mathcal G_i
        \right)
        \ge
        p_{\rm br}>0 .
\]
Consequently,
\[
        \mathbf P_{Q}(N>n)\le (1-p_{\rm br})^n,
        \qquad n\ge0 .
\]

Finally, the post-regeneration one-sided state is obtained by shifting
$(X_+,\tau_+)$ to $(0,0)$, retaining the distinguished leftmost barrier
$\Gamma_+$, and retaining the marked graphical construction on and to the
right of $\Gamma_+$.

\subsection{The barrier renewal property}
\label{subsec:barrier-renewal-property}

The key to renewal is that an accepted toy barrier
protects the fresh graphical environment from the
information revealed during the preceding exploration.
Consequently, the evolution retained beyond the barrier
has the same law at each accepted trial, independently
of the completed history.

We can now formulate the renewal statement. Its two ingredients are the
invariant trial law established in the preceding subsection and the
deterministic shielding property of an accepted toy path established in
Lemma~\ref{lem:toy-path-shielding}.

Since acceptance is determined by the infinite future, the renewal
statement is formulated with respect to the exploration $\sigma$-field
$\mathcal H_+$ defined below.

For each $i$, let $\mathcal K_i$ be as in
Lemma~\ref{lem:freshness-adaptive-trials}. On a successful preparation, write
\[
        S_i^{\rm toy}:=\{\mathcal E_i\in S^{\rm toy}\}.
\]
Recall the accepted-trial index $N$ from
\eqref{eq:one-cycle-first-accepted-index}.

Let $\mathcal H_+$ be the renewal $\sigma$-field characterized as follows:
on the event $\{N=i\}$ it contains all information in $\mathcal K_i$, the
fact that the $i$-th preparation succeeded, and the single bit of
information $S_i^{\rm toy}=1$, but it contains no further information from
$\mathcal E_i$. Equivalently, $B\in\mathcal H_+$ if, for each $i$, there
exists $B_i\in\mathcal K_i$ such that
\[
        B\cap\{N=i\}=B_i\cap\{N=i\}.
\]

\begin{proposition}[Barrier renewal property]
\label{prop:renewal}
The post-regeneration one-sided state $\Phi(\mathcal E_N)$ is independent of
$\mathcal H_+$ and has law $Q$. More precisely, for every bounded measurable
$G:\mathsf\Omega_+\to\mathbb R$,
\begin{equation}
\label{eq:renewal-conditional-law}
        \mathbf E_Q\!\left[
        G\bigl(\Phi(\mathcal E_N)\bigr)\mid\mathcal H_+
        \right]
        =
        Q(G)
        \qquad\mathbf{P}_{\mu}\text{ - a.s.}
\end{equation}
Moreover, the future CPTI right-edge process is determined by
$\Phi(\mathcal E_N)$ alone: any two admissible completions of the discarded
graphical history to the strict left of the accepted barrier produce the
same future right-edge process.
\end{proposition}
\begin{proof}
The distributional assertion is a consequence of the invariant trial-law lemma and
does not use the shielding argument.

Fix a bounded measurable $G$. Let $B\in\mathcal H_+$.  On $\{N=i\}$ write
\[
        B\cap\{N=i\}=B_i\cap\{N=i\},
\]
where $B_i$ is measurable with respect to the information available before
the accepted post-preparation environment is exposed.

By \eqref{eq:fresh-accepted-kernel}, on the event that the $i$-th
preparation succeeds,
\[
\mathbf E_{Q}\!\left[
        \Ind{S_i^{\rm toy}}
        G\bigl(\Phi(\mathcal E_i)\bigr)
        \,\middle|\,\mathcal K_i
\right]
=
\widehat\theta Q(G),
\]
while, taking $G\equiv1$,
\[
\mathbf P_{Q}\!\left(
        S_i^{\rm toy}\mid\mathcal K_i
\right)
=
\widehat\theta.
\]
Therefore
\begin{align}
&\mathbf E_{Q}\!\left[
  \Ind{B_i}\,
  \Ind{\{N>i-1\}}\,
  \Ind{\{\text{$i$-th preparation succeeds}\}}\,
  \Ind{S_i^{\rm toy}}\,
  G(\Phi(\mathcal E_i))
  \right]
\nonumber\\
&\qquad =
Q(G)\,
\mathbf E_{Q}\!\left[
  \Ind{B_i}\,
  \Ind{\{N>i-1\}}\,
  \Ind{\{\text{$i$-th preparation succeeds}\}}\,
  \Ind{S_i^{\rm toy}}
  \right].
\label{eq:renewal-factor-i}
\end{align}
The event multiplying $G$ on the left-hand side is precisely
$B\cap\{N=i\}$. Hence
\[
\mathbf E_Q\!\left[
    \Ind{B\cap\{N=i\}}\,
    G\bigl(\Phi(\mathcal E_i)\bigr)
\right]
=
Q(G)\,\mathbf P_Q\bigl(B\cap\{N=i\}\bigr).
\]
Summing over $i\ge1$ and using $N<\infty$ almost surely gives
\[\mathbf E_Q\!\left[
\Ind B\,G\bigl(\Phi(\mathcal E_N)\bigr)
\right]
=
Q(G)\,\mathbf P_Q(B).
\]
Since this holds for every $B\in\mathcal H_+$,
\eqref{eq:renewal-conditional-law} follows. Thus $\Phi(\mathcal E_N)$ has law $Q$
and is independent of $\mathcal H_+$.

It remains to verify that the retained one-sided state determines the future
right edge. Let two graphical histories in the strict-left region be
admissible completions of the same realized state $\Phi(\mathcal E_N)$. By definition,
the distinguished path contained in $\Phi(\mathcal E_N)$ is a legal CPTI infection
path after either completion is glued to the retained graphical construction.
The two resulting CPTI constructions have identical marked graphical
information on and to the right of this path. Therefore
Lemma~\ref{lem:toy-path-shielding} implies that their configurations agree
on and to the right of the barrier at every future time. In particular,
their right edges coincide. Hence the future right-edge process does not
depend on which admissible completion of the discarded strict-left history
is used, and is therefore determined by $\Phi(\mathcal E_N)$ alone.
\end{proof}

\subsection{Regeneration increments}
\label{subsec:regeneration-increments}

We now define the regeneration sequence.  Let $(X_1,\tau_1)$ be the first accepted barrier constructed from the original half-line initial condition.  In other words, $\tau_1$ is the time at which the first accepted preparation is completed, and $X_1$ is the prepared site.  After this time, the distinguished leftmost toy barrier launched from this site is used to define the post-regeneration one-sided state.

Inductively, after $(X_n,\tau_n)$ has been constructed, we shift this accepted barrier to $(0,0)$ and run the same record-level search from the one-sided state to the right of its distinguished leftmost toy barrier.  If the next accepted barrier, in shifted coordinates, has prepared site $X_+$ and completion time
$\tau_+$, we set
\[        X_{n+1}:=X_n+X_+,
        \qquad
        \tau_{n+1}:=\tau_n+\tau_+ .
\]
Thus $(X_n,\tau_n)$ is the space-time location of the $n$-th accepted barrier in the original coordinates.

Set $X_0:=0$ and $\tau_0:=0$.  For $n\ge1$, define the regeneration increments by
\[
        \Delta X_n:=X_n-X_{n-1},
        \qquad
        \Delta\tau_n:=\tau_n-\tau_{n-1}.
\]
We also define the maximal oscillation during the $n$-th cycle by
\[
        M_n:=
        \sup_{\tau_{n-1}\le t\le \tau_n}
        |r_t-X_{n-1}|.
\]
As a consequence of the barrier renewal property, the regeneration increments and maximal oscillations are i.i.d. after the first cycle; see Figure~\ref{fig:regeneration-increments}. This is made precise in the next proposition.

\begin{figure}
    \centering
    \scalebox{0.85}{\begin{tikzpicture}[x=0.85cm,y=0.85cm,>=Latex]

\tikzset{
    axis/.style={-{Latex[length=2.2mm]}, line width=0.7pt},
    guide/.style={densely dashed, gray!70, line width=0.6pt},
    edge/.style={
        black,
        thick,
        decorate,
        decoration={random steps, segment length=5pt, amplitude=0.55pt}
    },
    regpoint/.style={circle, fill=red!75!black, inner sep=1.5pt},
    lab/.style={font=\small},
    title/.style={font=\bfseries\large}
}


\fill[gray!3] (0.15,3.55) rectangle (10.2,5.85);

\fill[gray!6] (0.15,5.85) rectangle (10.2,8.45);


\def\MthreeLeft{3.45}
\def\MthreeRight{6.05}

\fill[gray!13] (\MthreeLeft,5.85) rectangle (\MthreeRight,8.45);

\draw[gray!55, dashed, line width=0.6pt]
    (\MthreeLeft,5.85) -- (\MthreeLeft,9.45);

\draw[gray!55, dashed, line width=0.6pt]
    (\MthreeRight,5.85) -- (\MthreeRight,9.45);

\draw[guide] (4.75,5.85) -- (4.75,9.45);

\draw[<->, thick, gray!75]
    (\MthreeLeft,9.10) -- (4.75,9.10);

\node[lab, above, gray!75!black] at (4.10,9.10) {$M_3$};


\draw[axis] (0,0) -- (0,10.2) node[above] {time};
\draw[axis] (0,0) -- (10.8,0) node[right] {space};


\draw[edge]
    (0.6,0.35)
    .. controls (0.95,0.95) and (1.30,1.35) .. (1.65,1.75)
    .. controls (1.30,2.15) and (1.85,2.55) .. (2.10,2.95)
    .. controls (2.45,3.25) and (2.35,3.45) .. (2.70,3.55)

    .. controls (2.55,3.95) and (2.35,4.20) .. (2.30,4.45)
    .. controls (2.65,4.90) and (3.35,5.20) .. (4.05,5.48)
    .. controls (4.40,5.62) and (4.65,5.76) .. (4.75,5.85)

    .. controls (4.55,6.15) and (3.95,6.45) .. (3.45,6.85)
    .. controls (3.85,7.25) and (4.65,7.55) .. (5.35,7.90)
    .. controls (5.65,8.08) and (5.85,8.30) .. (5.95,8.45)

    .. controls (6.25,8.90) and (7.10,9.20) .. (8.20,9.75);

\node[lab] at (1.25,2.25) {$r_t$};


\coordinate (R1) at (2.70,3.55);
\coordinate (R2) at (4.75,5.85);
\coordinate (R3) at (5.95,8.45);

\node[regpoint] at (R1) {};
\node[regpoint] at (R2) {};
\node[regpoint] at (R3) {};


\draw[guide] (0,3.55) -- (2.70,3.55);
\draw[guide] (2.70,0) -- (2.70,3.55);

\draw[guide] (0,5.85) -- (4.75,5.85);
\draw[guide] (4.75,0) -- (4.75,5.85);

\draw[guide] (0,8.45) -- (5.95,8.45);
\draw[guide] (5.95,0) -- (5.95,8.45);

\node[left]  at (0,3.55) {$\tau_1$};
\node[left]  at (0,5.85) {$\tau_2$};
\node[left]  at (0,8.45) {$\tau_3$};

\node[below] at (2.70,0) {$X_1$};
\node[below] at (4.75,0) {$X_2$};
\node[below] at (5.95,0) {$X_3$};

\node[red!75!black, lab, above left, xshift=-12pt, yshift=2pt]
    at (R1) {$(X_1,\tau_1)$};

\node[red!75!black, lab, above left, xshift=-4pt, yshift=2pt]
    at (4.63,5.67) {$(X_2,\tau_2)$};

\node[red!75!black, lab, above right, xshift=5pt, yshift=4pt, minimum width=37pt]
    at (6.45,8.25) {$(X_3,\tau_3)$};


\draw[<->, thick, gray!80]
    (0.55,3.55) -- (0.55,5.85);

\node[lab, right, gray!80!black] at (0.55,4.70)
{$\Delta\tau_2$};

\draw[<->, thick, red!75!black]
    (2.70,0.65) -- (4.75,0.65);

\node[red!75!black, lab, below] at (3.72,0.65)
{$\Delta X_2$};

\draw[<->, thick, red!75!black]
    (4.75,1.05) -- (5.95,1.05);

\node[red!75!black, lab, above] at (5.35,1.05)
{$\Delta X_3$};

\node[lab] at (1.66,0.51) {first cycle};
\node[lab] at (3.88,4.33) {cycle $2$};
\node[lab] at (5.35,7.18) {cycle $3$};

\node[
    lab,
    align=center,
    draw=red!75!black,
    rounded corners,
    fill=red!4,
    inner sep=3pt
] at (8.0,3.05)
{after the first cycle:\\
$(\Delta X_n,\Delta\tau_n,M_n)_{n\ge2}$ are i.i.d.};

\end{tikzpicture}}
    \caption{Regeneration cycles of the right edge. The second cycle illustrates the displacement
$\Delta X_2$ and duration $\Delta\tau_2$,
while the third illustrates the maximal
oscillation $M_3$.}
    \label{fig:regeneration-increments}
\end{figure}
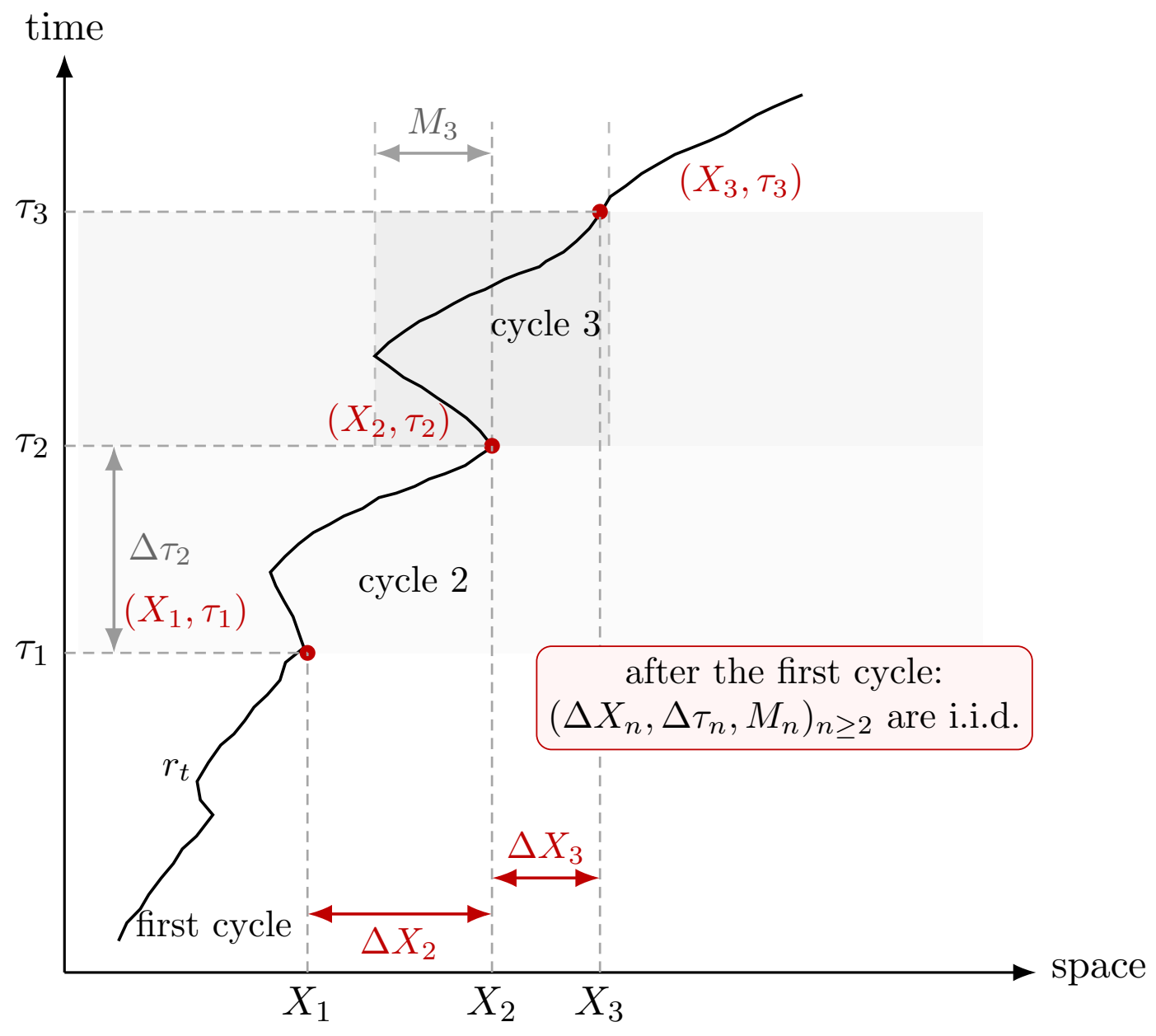
\begin{proposition} 
\label{prop:barrier-regeneration-increments}

After the first accepted barrier, the sequence
\begin{equation}
\label{eq:barrier-iid-increment-sequence}
        (\Delta X_n,\Delta\tau_n,M_n)_{n\ge2}
\end{equation}
is i.i.d. The first cycle may have a different law.
\end{proposition}

\begin{proof}
Once the first accepted barrier has been found, we translate its prepared site and time to $(0,0)$.  We then keep the marked graphical construction on and to the right of the distinguished leftmost toy barrier launched at the accepted preparation.  The same invariant trial-law argument as in Lemma~\ref{lem:freshness-adaptive-trials}, applied to the search from the original half-line initial condition, shows that at the first accepted trial the shifted one-sided state has law $Q$, independently of the information revealed before that accepted post-preparation environment is exposed.

Assume that after the $n$-th accepted barrier, with $n\ge1$, the shifted
one-sided state has law $Q$ and is independent of the completed past cycles. Run the same record-level search from this state.  The next cycle variables $(\Delta X_{n+1},\Delta\tau_{n+1},M_{n+1})$ are measurable functions of the one-sided state and of the graphical information revealed during the cycle.  Since the input law is $Q$, their conditional law, given the completed past, is the same fixed law for every $n\ge1$.  At the accepted barrier, Proposition~\ref{prop:renewal} gives a new independent post-regeneration one-sided state with law $Q$.  Iterating proves that the sequence in \eqref{eq:barrier-iid-increment-sequence} is i.i.d. 
\end{proof}

\section{Regenerative limit theory for the right edge}
\label{sec:lln-clt}

The renewal structure reduces the large-scale motion
of the right edge to a sequence of independent
regeneration cycles. To transfer limit theorems from
this regenerative skeleton to the full process, we
must also control the fluctuations of the right edge
between successive regeneration times.

The regenerative limit theorems rely on the following
exponential moment bounds, whose proof is deferred
to Appendix~\ref{app:moment-estimates}.
\begin{proposition}[Exponential moment bounds]
\label{prop:cycle-exponential-moments}
Fix the initial configuration in which every site $x \leq 0$ is infected and every site $x > 0$ is susceptible. There exists $\delta>0$ such that
\begin{equation}
\label{eq:first-cycle-moments}
\mathbf E_\mu\left[
e^{\delta\tau_1}
+e^{\delta|X_1|}
+e^{\delta M_1}
\right]<\infty.
\end{equation}
Moreover, for a generic regeneration cycle
$(X_+,\tau_+,M_+)$ started from $Q$,
\begin{equation}
\label{eq:cycle-exponential-moments}
\mathbf E_Q\left[
e^{\delta\tau_+}
+e^{\delta|X_+|}
+e^{\delta M_+}
\right]<\infty.
\end{equation}
\end{proposition}

To extend the regenerative limit theorems to arbitrary
$\xi\in\mathcal X^{\mathrm e}$, we first show that initial
immunities sufficiently far to the right expire before
the advancing infection can reach them.

\begin{lemma}[Entrance past the initial immunities]
\label{lem:entrance-initial-immunities}
Fix $\xi\in\mathcal X^{\mathrm{e}}$. For $y>0$, let $I_y$ be the residual
lifetime at time $0$ of the immunity at $y$ if $\xi(y)=-1$, and set
$I_y=0$ otherwise. Let
\[
\widehat T_y
:=
\inf\left\{
t\ge 0:
\text{there is a graphical arrow path from }
(-\infty,0]\times\{0\}
\text{ to }(y,t)
\right\}.
\]
Then, almost surely, there exists a finite random $R$ such that
\begin{equation}
\label{eq:initial-immunity-clears-before-arrival}
I_y<\widehat T_y
\qquad\text{for every }y\ge R.
\end{equation}
In particular, beyond the random level $R$, no immunity present at time
$0$ can affect the process before the corresponding site is first reached
by the infection.
\end{lemma}

\begin{proof}
Let $\widetilde T_y$ be the first arrival time at
$y$ when recovery marks are ignored. Then
$\widetilde T_y\le\widehat T_y$ and
$\widetilde T_y\stackrel d=E_1+\cdots+E_y$,
where the $E_i$ are i.i.d.\
$\operatorname{Exp}(\lambda)$. Hence, for some
$a,c,C>0$,
\[
\mathbf P_\mu(I_y\ge\widehat T_y)
\le Ce^{-cy}+e^{-\mu ay}.
\]
The bound is summable in $y$, so
Borel--Cantelli gives
\eqref{eq:initial-immunity-clears-before-arrival}.
\end{proof}

We now use the disappearance of distant initial immunities
to extend the regenerative construction to arbitrary
$\xi\in\mathcal X^{\mathrm e}$.

\begin{lemma}[Entrance into the regenerative regime]
\label{lem:entrance-regenerative-regime}
For every $\xi\in\mathcal X^{\mathrm e}$, the process
reaches an accepted barrier in finite time almost surely.
After translating this barrier to $(0,0)$, the retained
one-sided state has law $Q$, independently of the
information revealed before its post-preparation
environment is exposed.
\end{lemma}
\begin{proof}
Fix $\xi\in\mathcal X^{\mathrm e}$ and let
$A_\xi:=\{x:\xi(x)=1\}$, which is infinite.
By attractiveness of the toy process, positive survival from a singleton,
and a spatial restart argument using well-separated
sites of $A_\xi$, the toy process started from
$A_\xi$ survives almost surely. By
Remark~\ref{rem:toy-unbounded-growth}, its right
edge reaches arbitrarily large levels. Since every
toy infection path is a legal CPTI infection path,
the same holds for the CPTI started from $\xi$.

Run the adaptive trial construction of
Subsection~\ref{subsec:one-cycle-from-Q},
choosing each successive candidate beyond both the
running record and all previously revealed graphical
information, without reference to the random level $R$
of Lemma~\ref{lem:entrance-initial-immunities}.

Every rejected trial is resolved in finite time and
reveals information only up to a finite spatial level.
Since the right edge is unbounded, the candidate levels
tend to infinity. Almost surely, they therefore
eventually exceed $R$. By
Lemma~\ref{lem:entrance-initial-immunities}, inherited
immunities cannot interfere with the graphical
construction at these levels. The freshness argument
of Lemma~\ref{lem:freshness-adaptive-trials} then gives
the uniform conditional acceptance probability
$p_{\rm br}=q_{\rm prep}\widehat\theta>0$ for the
unaffected trials. Consequently, an accepted barrier
is reached in finite time almost surely.

Finally, Proposition~\ref{prop:renewal} shows that selection
of the first accepted trial does not bias its
post-preparation state. After translation, this
state has law $Q$ and is independent of the
information revealed before its post-preparation
environment is exposed.
\end{proof}

We now derive the limit theorems from the regeneration
structure.

\begin{proof}[Proof of Theorem~\ref{teo_lln}]
By Propositions~\ref{prop:barrier-regeneration-increments}
and~\ref{prop:cycle-exponential-moments}, the post-first
cycles are i.i.d.\ and integrable. The regenerative
law of large numbers
\cite[Theorem~VI.3.1]{Asmussen2003} therefore gives
\[
\frac{r_t}{t}\longrightarrow
v:=\frac{\mathbf E_Q[X_+]}{\mathbf E_Q[\tau_+]}
\qquad \mathbf P_Q\text{-a.s.}
\]
Since each cycle reaches a new record level,
$X_+>0$ almost surely, $v>0$.

By Lemma~\ref{lem:entrance-regenerative-regime},
every process started from $\xi\in\mathcal X^{\mathrm e}$
enters the regenerative regime at an almost surely finite
time. Thereafter, the regeneration cycles have the same
law as under $Q$. Since the initial time and spatial
displacement are finite almost surely, they do not affect
the limiting speed. Thus the law of large numbers holds
for every such $\xi$.

For the half-line initial condition,
\eqref{eq:first-cycle-moments} and
\eqref{eq:cycle-exponential-moments} give uniform
integrability of $r_t/t$. Thus the convergence
also holds in $L^1$.
\end{proof}

\begin{proof}[Proof of Proposition~\ref{prop:speed-limit-mu-infty}]
Couple the CPTI and classical contact process, both
started from the half-line, using the same infection
arrows and recovery marks. Since immunity only
suppresses infections, $r_t^\mu\le r_t^{\rm cl}$, with $r_t^{\rm cl}$ being the right edge of the classical contact process.
By~\eqref{eq:right-edge-lln} and the classical
right-edge law of large numbers,
\begin{equation}
\label{eq:speed_upper_bound}
\limsup_{\mu\to\infty}v(\lambda,\mu)
\le v_{\rm CP}(\lambda).
\end{equation}

For the reverse inequality, fix $\varepsilon>0$ and
choose the classical block event $G_{j,n}^{\rm cl}$
as in Remark~\ref{rem:classical-speed-blocks}.
Choose $p_\star<1$ sufficiently close to one that
the corresponding oriented percolation speed, in
microscopic coordinates, is at least
$v_{\rm CP}(\lambda)-\varepsilon$.

Choose $\rho_\star>0$ so that
Proposition~\ref{prop:peierls-after-old-tail-deletion}
applies whenever $p\ge p_\star$ and
$\rho_E\le\rho_\star$. By the cutoff approximation
and the old-tail estimate of
Subsections~\ref{subsec:cutoff-block-construction}
and~\ref{subsec:old-tail-obstacles}, we may choose
the cutoff and then $\mu$ sufficiently large that
both conditions hold.

Proposition~\ref{prop:peierls-after-old-tail-deletion}
and Remark~\ref{rem:classical-speed-blocks} then give,
with positive probability,
\begin{equation} 
\label{eq:lim-inf-speed}
\liminf_{t\to\infty}\frac{r_t^\mu}{t}
\ge v_{\rm CP}(\lambda)-\varepsilon.
\end{equation}
Using Theorem~\ref{teo_lln} and~\eqref{eq:lim-inf-speed},
\[
\liminf_{\mu\to\infty}v(\lambda,\mu)
\ge v_{\rm CP}(\lambda).
\]
Together with~\eqref{eq:speed_upper_bound}, this
proves the proposition.
\end{proof}

\begin{proof}[Proof of Theorem~\ref{teo_fclt}]
By Proposition~\ref{prop:barrier-regeneration-increments},
the post-first cycles are i.i.d. Moreover,
\eqref{eq:cycle-exponential-moments} gives
\[
u^2\mathbf P_Q(M_++v\tau_+>u)\longrightarrow0
\qquad\text{as }u\to\infty.
\]
Thus the regenerative functional central limit theorem
\cite[Theorem~1.(b)]{GlynnWhitt1993} applies.

The limiting process is Brownian motion with variance
parameter
\[
\sigma^2
=
\frac{\mathbf E_Q[(X_+-v\tau_+)^2]}
{\mathbf E_Q[\tau_+]}.
\]
By Lemma~\ref{lem:entrance-regenerative-regime},
every process started from $\xi\in\mathcal X^{\mathrm e}$
enters the regenerative regime at an almost surely finite
time. After this entrance, the regeneration cycles have
the same law as under $Q$. The finite initial time and
spatial displacement vanish under diffusive scaling,
so the functional limit holds for every such $\xi$.

Finally, $\sigma^2>0$. Indeed, on a positive-probability
finite graphical event, $X_+$ is fixed while a Poisson
waiting time contributing to $\tau_+$ varies over a
nondegenerate interval without changing the regeneration
construction. Since $v>0$, $X_+-v\tau_+$ is not almost
surely constant, proving the claim.
\end{proof}

\section{Shielding and asymptotic laws}
\label{sec:shielding-asymptotic-laws}

We now use shielding to couple CPTI evolutions with
different initial configurations and graphical histories.
The key idea is that common toy barriers force agreement
in expanding space-time regions, despite the lack of
attractiveness. Combined with the regenerative structure,
this allows us to study the process seen from the right
edge and establish complete convergence.

\subsection{Regenerative coupling and shielding}
\label{subsec:regenerative-coupling-shielding} 
From now on, we will consider not only configurations with infinitely many infected individuals to the left of the origin, but also configurations with infinitely many infected individuals on both sides of the origin.
Once the two processes have entered the regenerative regime,
their regeneration structures can be successfully coupled.
Shielding then allows the common graphical environment
to determine their evolution in expanding regions, despite
differences in their previous histories. For convergence
seen from the right edge, the regeneration sites are first
translated to the origin. For bulk convergence, the coupling
is performed in the original spatial coordinates, using
barriers expanding in both directions. 

\begin{lemma}[One- and two-sided regenerative coupling]
	\label{lem:finite-window-coupling}
	For $i=1,2$, let $(\xi_t^{(i)})_{t\ge0}$ be a CPTI whose first
	accepted regeneration time is almost surely finite, and let
	$(X_n^{(i)},\tau_n^{(i)})_{n\ge1}$ be its regeneration sequence.
	
	\begin{enumerate}
		\renewcommand{\labelenumi}{(\roman{enumi})}
		\renewcommand{\theenumi}{(\roman{enumi})}
		
		\item\label{item:coupling-one-sided}
		The two processes may be coupled so that there exist an almost surely
		finite common regeneration epoch $\tau_{\mathrm{cpl}}$ and corresponding
		regeneration sites $X_{\mathrm{cpl}}^{(1)}$ and
		$X_{\mathrm{cpl}}^{(2)}$. After translating these sites to the origin,
		the two post-regeneration one-sided graphical states are identified and
		the same future graphical construction is used in the retained region.
		For $s\ge0$, set
		\[
		\tilde{\xi}_s^{(i)}(x)
		:=
		\xi_{\tau_{\mathrm{cpl}}+s}^{(i)}
		\bigl(X_{\mathrm{cpl}}^{(i)}+x\bigr),
		\qquad i=1,2.
		\]
		
		Let $\Gamma^-$ be the distinguished leftmost infinite toy infection
		path of the common post-regeneration one-sided graphical state. Then
		\begin{equation}
			\label{eq:coupling-recentered}
			\tilde{\xi}_s^{(1)}(x)
			=
			\tilde{\xi}_s^{(2)}(x),
			\qquad x\ge\Gamma^-(s),\quad s\ge0.
		\end{equation}
		Moreover, for some $c_->0$,
		\[
		\limsup_{s\to\infty}
		\frac{\Gamma^-(s)}{s}
		\le -c_-
		\qquad\text{almost surely}.
		\]
		
		Setting $D:=X_{\mathrm{cpl}}^{(2)}-X_{\mathrm{cpl}}^{(1)}$, we have
		\begin{equation}
			\label{eq:coupling-with-displacement}
			\xi_{\tau_{\mathrm{cpl}}+s}^{(1)}(x)
			=
			\xi_{\tau_{\mathrm{cpl}}+s}^{(2)}(x+D),
			\qquad
			x\ge X_{\mathrm{cpl}}^{(1)}+\Gamma^-(s).
		\end{equation}
		
		Furthermore, for every $\ell<\infty$, there exists an almost surely
		finite $T_\ell$ such that
		\[
		\Xi_{\tau_{\mathrm{cpl}}+s}^{(1)}(x)
		=
		\Xi_{\tau_{\mathrm{cpl}}+s}^{(2)}(x),
		\qquad -\ell\le x\le0,
		\]
		for all $s\ge T_\ell$.
		
		\item\label{item:coupling-two-sided}
		The two processes may instead be coupled so that there exist an almost
		surely finite coupling epoch $\tau_{\mathrm{cpl}}$ and a site
		$X_{\mathrm{cpl}}$ such that, at
		$(X_{\mathrm{cpl}},\tau_{\mathrm{cpl}})$, their post-regeneration
		one-sided graphical states are identified in the same spatial
		coordinates and the same future graphical construction is used.
		
		Let $\Gamma^-$ and $\Gamma^+$ be respectively the leftmost and
		rightmost infinite toy infection paths issued from the prepared site.
		For some $c_-,c_+>0$,
		\[
		\limsup_{s\to\infty}\frac{\Gamma^-(s)}{s}\le-c_-,
		\qquad
		\liminf_{s\to\infty}\frac{\Gamma^+(s)}{s}\ge c_+
		\qquad\text{almost surely}.
		\]
		Consequently,
		\[
		\liminf_{s\to\infty}
		\frac{\Gamma^+(s)-\Gamma^-(s)}{s}
		\ge c_-+c_+>0
		\qquad\text{almost surely}.
		\]
		
		The strict-left side of $\Gamma^-$ may be replaced by any admissible
		completion, and the strict-right side of $\Gamma^+$ may be replaced by
		any admissible completion in the spatially reflected sense, separately
		or simultaneously. Under any such replacements, the resulting CPTI
		evolution agrees with the original one whenever
		\[
		X_{\mathrm{cpl}}+\Gamma^-(s)
		\le x\le
		X_{\mathrm{cpl}}+\Gamma^+(s).
		\]
		
		In particular, for every $\ell<\infty$, there exists an almost surely
		finite $T_\ell$ such that
		\[
		\xi_{\tau_{\mathrm{cpl}}+s}^{(1)}(x)
		=
		\xi_{\tau_{\mathrm{cpl}}+s}^{(2)}(x),
		\qquad -\ell\le x\le\ell,
		\]
		for all $s\ge T_\ell$.
	\end{enumerate}
\end{lemma}

\begin{figure}[t]
	\input{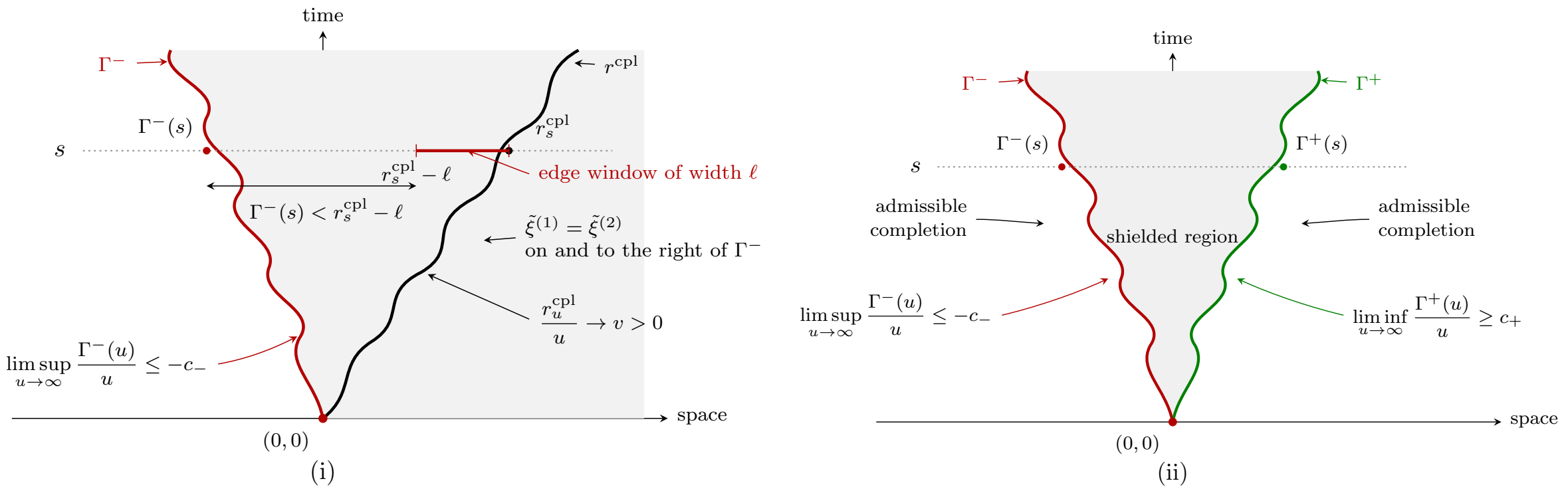}
	\caption{Regenerative couplings of
Lemma~\ref{lem:finite-window-coupling}.
(i) One-sided coupling after translating the regeneration
sites to the origin, with agreement to the right of
$\Gamma^-$.
(ii) Two-sided coupling from a common regeneration site,
with agreement between $\Gamma^-$ and $\Gamma^+$.}
	\label{fig:coupling}
\end{figure}

\begin{proof}
	Figure~\ref{fig:coupling} illustrates the two coupling mechanisms used below.
	
We first prove~\ref{item:coupling-one-sided}.
After their first accepted regeneration times, the
two processes have the same regenerative law, with
finite mean cycle length by
Proposition~\ref{prop:cycle-exponential-moments}.
Moreover, the inter-regeneration time has a spread-out
distribution: on a positive-probability graphical event,
one of the Poisson waiting times can vary over a
nondegenerate interval without changing the regeneration
construction. Hence, by
\cite[Theorem~1.1]{Thorisson1983}, the regeneration
sequences can be coupled, after factoring out spatial
translations, to have an almost surely finite common
regeneration epoch $\tau_{\mathrm{cpl}}$. Let $X_{\mathrm{cpl}}^{(1)}$ and
$X_{\mathrm{cpl}}^{(2)}$ be the corresponding regeneration sites.
By Proposition~\ref{prop:renewal}, after translating these sites to
the origin, the two post-regeneration one-sided graphical states have
law $Q$, independently of their respective pasts. We identify these
states and use the same future graphical marks in the retained region.
	
Let $\Gamma^-$ be the distinguished leftmost infinite toy infection
path of this common post-regeneration state. By the shielding
property of Lemma~\ref{lem:toy-path-shielding},
	\[
	\tilde{\xi}_s^{(1)}(x)
	=
	\tilde{\xi}_s^{(2)}(x),
	\qquad x\ge\Gamma^-(s),\quad s\ge0.
	\]
This proves~\eqref{eq:coupling-recentered}, and returning to the
original coordinates gives~\eqref{eq:coupling-with-displacement}.
	
By Proposition~\ref{prop:peierls-after-old-tail-deletion}
and Remark~\ref{rem:toy-unbounded-growth},
$\Gamma^-$ has strictly negative asymptotic speed.
Since the two shifted configurations agree to the
right of $\Gamma^-$, their right edges coincide;
denote the common edge by $r_s^{\mathrm{cpl}}$.
By Theorem~\ref{teo_lln}, this edge has asymptotic
speed $v>0$. Consequently,
\[
    r_s^{\mathrm{cpl}}-\Gamma^-(s)
    \longrightarrow\infty
    \qquad\text{almost surely}.
\]
Thus, every fixed window behind the common right edge
eventually lies in the region of agreement, proving
\ref{item:coupling-one-sided}.
	
We now prove~\ref{item:coupling-two-sided}.
The successive regeneration points form a random walk on
$\mathbb Z\times\mathbb R$ with i.i.d. increments
$(\Delta X,\Delta\tau)$. The graphical construction shows
that every positive spatial increment is possible and that,
for consecutive spatial increments, the corresponding
regeneration-time distributions have overlapping absolutely
continuous components. These overlaps, together with
continuous variation of the regeneration times, give the
joint aperiodicity and spread-outness needed for the
random-walk coupling argument of
\cite[Proposition~1]{LindvallRogers1996}. That argument,
applied on $\mathbb Z\times\mathbb R$, couples the
regeneration sequences at an almost surely finite common
space-time point $(X_{\mathrm{cpl}},\tau_{\mathrm{cpl}})$.

Since the rightmost infinite toy path lies to the right
of the distinguished leftmost path, both extremal paths
are determined by the common retained graphical state.
Thus the same pair of shielding paths is available
to both coupled evolutions.

At the common regeneration point, the post-regeneration graphical states can be identified, since both have law $Q$ independently of their respective completed pasts. The subsequent marked regenerative evolutions can therefore be taken identical.

	Let $\Gamma^-$ and $\Gamma^+$ be the leftmost and
rightmost infinite toy infection paths issued from
the prepared site. The leftward speed estimate above
and its spatially reflected counterpart give the
asserted linear expansion of the interval between
these paths.

By Lemma~\ref{lem:toy-path-shielding}, replacing the
strict-left graphical history of $\Gamma^-$ by an
admissible completion leaves the CPTI evolution
unchanged on and to its right. The reflected lemma
gives the analogous property to the left of
$\Gamma^+$. These replacements may be performed
successively, since each preserves the evolution
between the two paths. Thus, identifying the
exterior histories of the two coupled processes as
admissible completions yields agreement whenever
\[
X_{\mathrm{cpl}}+\Gamma^-(s)
\le x\le
X_{\mathrm{cpl}}+\Gamma^+(s).
\]
Since the barriers expand linearly in opposite
directions, this interval eventually contains every
fixed finite interval, proving
\ref{item:coupling-two-sided}.
\end{proof}

\subsection{The invariant edge law}
\label{subsec:invariant-edge-law}

The renewal structure allows us to construct a stationary
sequence of marked regeneration cycles. By extending this
construction backward in time, shielding ensures that the
configuration in each fixed window behind the right edge
eventually becomes independent of how far in the past the
process was started. This yields a stationary edge
configuration, while the same coupling mechanism gives
convergence from arbitrary initial configurations.

\begin{proof}[Proof of Theorem~\ref{teo:edge-invariant-law}]
By Proposition~\ref{prop:cycle-exponential-moments},
the regeneration times have finite mean increments.
Together with Proposition~\ref{prop:renewal}, this
allows us to construct a bi-infinite stationary
version of the marked regenerative process, whose
successive cycles have the prescribed regenerative
law. By abuse of notation, denote its law by $\mathbf P_Q$.
	
For $n\ge1$, let $(X_{-n},\tau_{-n})$ denote the $n$-th regeneration
point before time $0$, and let $\Gamma^{(-n)}$ be the distinguished
barrier issued from this point. Using the common marked graphical
construction, concatenate the successive regenerative pieces so that
each barrier is preserved by all constructions started at earlier
regeneration points. Let $\Xi^{(-n)}$ be the process seen from its
right edge obtained by starting from the regenerative state at
$(X_{-n},\tau_{-n})$ and evolving forward to time $0$.

By Lemma~\ref{lem:finite-window-coupling}
\ref{item:coupling-one-sided}, for every $\ell\ge1$,
\[
\mathbf P_Q\bigl(
r_0-\Gamma^{(-n)}(0)\le\ell
\bigr)\longrightarrow0.
\]
We may therefore choose a subsequence $(n_k)$ such that
the preceding probability, with $n=n_k$ and $\ell=k$,
is at most $2^{-k}$. By Borel--Cantelli,
\[
r_0-\Gamma^{(-n_k)}(0)>k
\]
for all sufficiently large $k$, almost surely.

Fix $\ell\ge1$. For all sufficiently large $k$,
the window $\{r_0-\ell,\ldots,r_0\}$ lies on or to
the right of $\Gamma^{(-n_k)}(0)$. By the shielding
property of Lemma~\ref{lem:toy-path-shielding},
the restrictions
\[
\Xi_0^{(-n_k)}\big|_{\{-\ell,\ldots,0\}}
\]
eventually agree. Since these restrictions are
compatible as $\ell$ varies, they define a random
edge configuration $\Xi_0^\star$.

Moreover, $\Xi_0^\star$ has infected sites arbitrarily far behind
its edge. Indeed, fix $j$. By the concatenated construction and
Lemma~\ref{lem:toy-path-shielding}, for every $k\ge j$,
\[
\Xi_0^{(-n_k)}
\bigl(\Gamma^{(-n_j)}(0)-r_0\bigr)=1.
\]
Passing to the limit with $j$ fixed gives
\[
\Xi_0^\star
\bigl(\Gamma^{(-n_j)}(0)-r_0\bigr)=1.
\]
Since $r_0-\Gamma^{(-n_j)}(0)\to\infty$ along the
subsequence, these infected sites extend arbitrarily far
to the left. Hence
$\Xi_0^\star\in\mathcal X^{\mathrm e}$ a.s.

	Let $\nu^{\mathrm e}$ denote the law of
$\Xi_0^\star$. By stationarity of the marked
regenerative construction, the same backward
procedure defines $\Xi_t^\star$ for every
$t\in\mathbb R$, with
$\mathcal L(\Xi_t^\star)=\nu^{\mathrm e}$.
		
	We next prove convergence. Fix
$\xi\in\mathcal X^{\mathrm e}$. By
Lemma~\ref{lem:entrance-regenerative-regime},
the process started from $\xi$ enters the regenerative
regime in finite time almost surely. Since
$(\Xi_t^\star)_{t\ge0}$ is constructed from the
stationary marked regeneration process with the same
regenerative law, the coupling of
Lemma~\ref{lem:finite-window-coupling}
\ref{item:coupling-one-sided} applies to these two
processes. Consequently, for every $\ell<\infty$,
there is an almost surely finite $T_\ell$ such that
\[
\Xi_t(x)=\Xi_t^\star(x),
\qquad -\ell\le x\le0,
\]
for all $t\ge T_\ell$.	
	If $f$ is a bounded local function depending only on
	$\{-\ell,\ldots,0\}$, then $f(\Xi_t)=f(\Xi_t^\star)$ for all
	sufficiently large $t$, almost surely. Since
	$\Xi_t^\star\sim\nu^{\mathrm{e}}$ for every $t$, bounded convergence
	gives
	\[
	\lim_{t\to\infty}
	\mathbb E^\xi_{\lambda,\mu}[f(\Xi_t)]
	=
	\int f\,d\nu^{\mathrm{e}}.
	\]
    Hence, the weak convergence in~\eqref{eq:weak-convergence-edge} holds.
    
	Since $(\Xi_t^\star)_{t\in\mathbb R}$ is a stationary
version of the edge process, its marginal law
$\nu^{\mathrm e}$ is invariant. Uniqueness follows
from convergence to $\nu^{\mathrm e}$ from every
initial configuration in $\mathcal X^{\mathrm e}$
and dominated convergence.
\end{proof}

\begin{remark}[Uniform convergence from finite initial sets]
\label{rem:uniform-edge-convergence-strategy}
We briefly indicate how the proof of
Theorem~\ref{teo:edge-invariant-law} could be adapted
to obtain the uniform convergence stated in
Remark~\ref{rem:uniform-edge-convergence}.

The main difference is that the initial infected set
is now finite and convergence is required uniformly
over sets $A$ with $\operatorname{diam}(A)\le ct$,
conditional on survival up to time $t$.
The first step would be to show that, under this
conditioning, the process enters the regenerative
regime sufficiently early, uniformly over such sets.
This could be achieved using the regeneration trials
of Section~\ref{sec:regeneration-renewal} and the corresponding quantitative estimates. Once regeneration has occurred, the coupling
of Lemma~\ref{lem:finite-window-coupling}%
~\ref{item:coupling-one-sided}, together with the
renewal estimates, would provide uniform control
of the time needed for agreement with the stationary
edge process in any fixed window.
The condition $c<v$ would enter in controlling
the initial transient relative to the linear growth
of the infection front.
The shielding argument used in the proof of
Theorem~\ref{teo:edge-invariant-law} would then
yield the desired convergence.
We omit the details.
\end{remark}

\subsection{Complete convergence}
\label{subsec:complete-convergence}
We now turn to complete convergence. Unlike convergence
seen from the edge, this requires controlling the
process in fixed spatial windows, as well as the
possibility of extinction. The proof combines
entrance into the regenerative regime with shielding
couplings from both directions.

\begin{lemma}[Entrance into the regenerative regime]
\label{lem:finite-entrance-regeneration}
Let $\xi\in\{-1,0,1\}^{\mathbb Z}$ have finite
nonempty infected set and finitely many immune sites,
and set $S_\xi=\{\tau^\xi=\infty\}$.
On $S_\xi$, the process enters the regenerative
regime in finite time almost surely, to the right
or to the left.

After translating the accepted regeneration point
to $(0,0)$ and reflecting space if necessary,
the retained one-sided graphical state has law $Q$,
independently of the previously revealed information.
\end{lemma}

\begin{proof}
On $S_\xi$, write
\[
    \ell_t:=\inf I_t^\xi,
    \qquad
    r_t:=\sup I_t^\xi.
\]
We first observe that survival forces the infection
to leave every finite interval. Indeed, the CPTI
restricted to $[-m,m]$, with infection attempts
leaving the interval suppressed, is a finite-state
Markov chain that almost surely reaches a
configuration without infections. Since the restricted
and unrestricted processes agree until the first
infection outside $[-m,m]$,
\[
\mathbb P_{\lambda,\mu}
\bigl(
S_\xi,\,
I_t^\xi\subset[-m,m]\text{ for all }t\ge0
\bigr)=0.
\]
Taking the union over $m$ gives
\begin{equation}
\mathbb P_{\lambda,\mu}\left(
S_\xi\cap
\left\{
\sup_{t\ge0}r_t<\infty,\;
\inf_{t\ge0}\ell_t>-\infty
\right\}
\right)=0.
\label{eq:finite-survival-unbounded}
\end{equation}	
Suppose first that $\sup_{t\ge0}r_t=\infty$.
We apply the adaptive record-level trial construction
of Subsection~\ref{subsec:one-cycle-from-Q}.
After each rejected trial, the next candidate level
is chosen beyond the running record and all previously
revealed graphical information. By
Lemma~\ref{lem:freshness-adaptive-trials},
each trial is therefore accepted, conditionally
on the past, with probability at least
\[
p_{\rm br}=q_{\rm prep}\widehat\theta>0.
\]
Since every rejected trial is resolved in finite
time almost surely, the probability that the first
$n$ trials are all rejected is at most
$(1-p_{\rm br})^n$. Hence a regeneration occurs
almost surely on the event
$\{\sup_{t\ge0}r_t=\infty\}$.

The reflected construction applies when the infection
escapes to the left. Together with~\eqref{eq:finite-survival-unbounded}, this proves the entrance assertion. The statement concerning
the retained graphical law follows from
Proposition~\ref{prop:renewal}, with reflection
in the left-hand case.
\end{proof}

We prove complete convergence by comparing finite-start
processes with the nontrivial invariant law of~\cite{AndjelSchinazi1996}. Finite approximations, two-sided shielding,
and the regenerative coupling identify this law as the
limit conditional on survival.

\begin{proof}[Proof of Theorem~\ref{teo:complete-convergence}]
Let $\nu$ be the nontrivial translation-invariant
invariant measure of Andjel and Schinazi
\cite{AndjelSchinazi1996}, normalized to have no
absorbing component. By translation invariance,
a configuration sampled from $\nu$ cannot have
a finite nonempty infected set. Stationarity and
the absence of an absorbing component then imply
that it has infinitely many infected sites almost
surely.

Let $\zeta_0$ have distribution $\nu$, and let
$\zeta^{(n)}_0$ agree with $\zeta_0$ on $[-n,n]$
and be susceptible elsewhere. Couple their
evolutions using the same graphical construction.
Since $\zeta_0$ has infinitely many infected sites almost
surely, we may condition on $\zeta_0$ and select arbitrarily
many well-separated infected sites. The toy graphical
construction is independent of $\zeta_0$. Positive toy
survival, reflection symmetry, and spatial mixing therefore
imply that, with probability tending to one as $n\to\infty$,
there are toy infection paths starting from infected sites
in $[-n,n]$ and expanding in opposite directions. These paths may be chosen with the leftward-expanding
path starting to the left of the origin and the
rightward-expanding path starting to its right, so
that their starting points enclose the origin.
 Since the initial configurations
agree between these paths, Lemma~
\ref{lem:toy-path-shielding} gives, for every
finite $\Lambda\subset\mathbb Z$,
\[
\lim_{n\to\infty}\limsup_{t\to\infty}
\mathbb P\bigl(
\zeta_t^{(n)}|_\Lambda
\ne\zeta_t|_\Lambda
\bigr)=0.
\]

In fact, these paths are common legal infection paths for both
evolutions. Since the initial configurations agree
between their starting points, Lemma~\ref{lem:toy-path-shielding}, applied
to each path (with spatial reflection for the
right-hand barrier), shows that the two evolutions
agree at all subsequent times in the region between
the paths. By ballisticity, this region eventually
contains every fixed finite interval.
In particular, the survival probability of
$\zeta^{(n)}$ tends to one.

Now fix finite nonempty $A\subset\mathbb Z$ and
write $S_A=\{\tau^A=\infty\}$. The proof of
Lemma~\ref{lem:finite-entrance-regeneration}
applies equally to finite-support initial
configurations with immune sites. Thus, for
each fixed $n$, Lemma~
\ref{lem:finite-window-coupling}
\ref{item:coupling-two-sided} couples the
surviving processes started from $A$ and
$\zeta^{(n)}_0$ so that they eventually agree
on every fixed finite interval.

Since $\zeta_t$ has distribution $\nu$ for every
$t$, the preceding approximation and coupling
imply, by letting first $t\to\infty$ and then
$n\to\infty$, that
\[
\mathbb P_{\lambda,\mu}
\bigl(\xi_t^A\in\cdot\mid S_A\bigr)
\Longrightarrow\nu.
\]
On $S_A^c$, the infection dies out in finite time
and the remaining immunity intervals expire.
Consequently,
\[
\mathbb P_{\lambda,\mu}(\xi_t^A\in\cdot)
\Longrightarrow
\mathbb P_{\lambda,\mu}(S_A^c)\delta_{\mathbf0}
+\mathbb P_{\lambda,\mu}(S_A)\nu.
\]

Finally, let $A$ be infinite and choose finite
sets $A_n\uparrow A$. By
Proposition~\ref{prop:peierls-after-old-tail-deletion},
reflection symmetry, and FKG, together with
spatial mixing, the probability that two toy
infection paths expanding in opposite directions
arise from sites of $A_n$ tends to one.
Coupling the processes from $A_n$ and $A$
through the same graphical construction,
Lemma~\ref{lem:toy-path-shielding} gives
eventual agreement on every fixed finite
interval, outside an event of probability
tending to zero as $n\to\infty$.

These paths also imply that
$\mathbb P_{\lambda,\mu}(S_A)=1$ and
$\mathbb P_{\lambda,\mu}(S_{A_n})\to1$.
Complete convergence for finite initial sets,
followed by $n\to\infty$, therefore yields
\[
\mathbb P_{\lambda,\mu}(\xi_t^A\in\cdot)
\Longrightarrow\nu.
\]
\end{proof}

\section{Extinction in the subcritical return regime}
\label{sec:extinction}

The extinction argument relies on the observation that
survival from a finite infected set requires repeated
infections of previously visited sites. The key is to
control these reinfections: each infection of a new site
creates a possible route for the infection to return,
but such a return requires the former source to recover,
lose its immunity, and be reinfected before its neighbor
recovers. This mechanism is reminiscent of the
Markov-modulated point-process comparison of
Broman~\cite{Broman2007}, although our argument uses
conditional first-moment estimates rather than
stochastic domination.

\begin{proof}[Proof of Theorem~\ref{teo_dies}]
	An infection event at a site $x$ is called a
\emph{return infection} if $x$ was infected at some
earlier time, including time zero. Otherwise, it is
called a \emph{first infection}. By~\cite[Theorem~5]{Heeszel2026}, survival from a finite infected set implies infinitely many
return infections. Thus, it suffices to show
that the number of return infections is almost
surely finite.
	
	We organize the return infections into generations according to their
	infection genealogy. A return infection belongs to generation $1$ if its
	ancestral line contains no earlier return infection. Inductively, a return
	infection belongs to generation $m+1$ if the most recent return infection on
	its ancestral line belongs to generation $m$. Let $\mathcal R_m$ denote the
	set of return infections in generation $m$, and write $Z_m:=|\mathcal R_m|$.
	
	Starting from the space-time points corresponding to $\mathcal R_m$, explore
	all legal descendants, stopping each branch when the next return infection is
	encountered. Let $\mathcal A_m$ be the collection of first-infection
	edges discovered during this exploration, and $G_m$ the number of return gates at generation $m$.
	
We first count the possible gates leading to
$\mathcal R_{m+1}$. Consider an infection edge
$u\to v$ discovered during the exploration.
Immediately after this infection, both sites are
infected. A return infection may subsequently occur
either at $u$, through an arrow from $v$, or at $v$,
through an arrow from $u$. In either case, the site
receiving the infection must first recover and lose
its immunity, while the other site must remain
infected until the return occurs.

Thus, each infection edge gives rise to two possible
directed return gates. Counting also the edges
corresponding to $\mathcal R_m$, we obtain
\begin{equation}
    G_m\leq 2(Z_m+|\mathcal A_m|).
    \label{eq:number-gates}
\end{equation}
	
	Fix an element of $\mathcal R_m$ and one spatial
direction. Each successive first infection requires
an infection arrow before the next recovery mark,
which occurs with probability at most
$\lambda/(1+\lambda)$. Thus, the expected number
of first-infection edges in that direction is at most
\[
\sum_{n\geq1}
\left(\frac{\lambda}{1+\lambda}\right)^n
=\lambda.
\]
Summing over both directions and all elements of
$\mathcal R_m$ yields
\begin{equation}
\mathbf E_{\lambda,\mu}|\mathcal A_m|
\leq 2\lambda\,\mathbf E_{\lambda,\mu}Z_m.
\label{eq:first-infection-count}
\end{equation}
	
	It remains to estimate the probability that a candidate gate succeeds.
	Consider a candidate gate associated with an
infection edge at time $s$. Let $u$ denote the
potential recipient of the return infection and
$v$ its potential source. Immediately after the
edge is created, both sites are infected.
	We evaluate the gate at its creation time. Conditional on the graphical
	history up to time $s$, let $R_u,R_v\stackrel{\mathrm{iid}}{\sim}\operatorname{Exp}(1)$
	be the residual recovery times of $u$ and $v$, let $I_u\sim\operatorname{Exp}(\mu)$
	be the immunity time of $u$ after its recovery, and let	$A_{v,u}\sim\operatorname{Exp}(\lambda)$
	be the waiting time, after the immunity of $u$ is lost, for the next infection
	arrow $v\to u$. These random variables are independent.
	
	For this gate to produce the next return infection, it is necessary that
	\[
	R_u+I_u+A_{v,u}<R_v.
	\]
	The geometry of the exploration and this reverse-gate mechanism are
	illustrated in Figure~\ref{fig:return-gate-exploration}.
	
	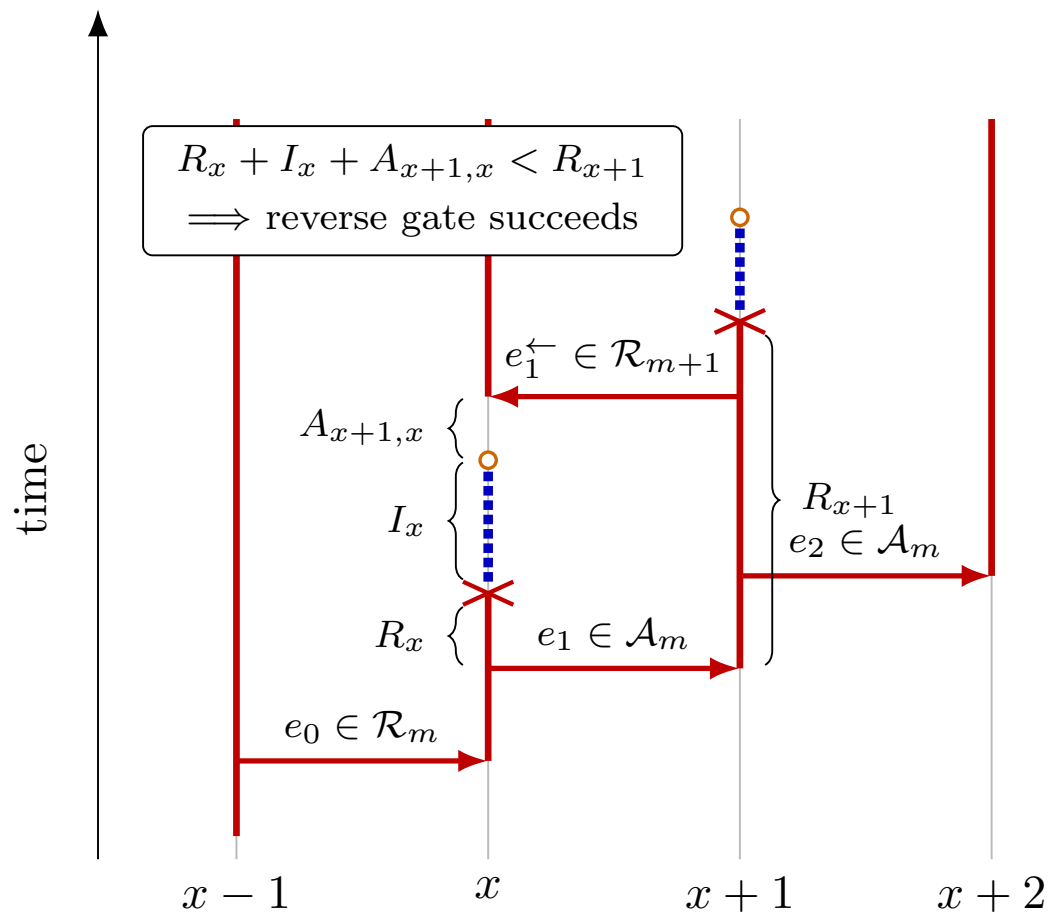
\begin{figure}[t]
		\centering
		\scalebox{0.9}{\centering
\begin{tikzpicture}[
    x=1.85cm,
    y=0.85cm,
    site/.style={
        gray!55,
        line width=0.4pt
    },
    infection/.style={
        -{Latex[length=2.2mm]},
        red!75!black,
        line width=1.1pt
    },
    infected/.style={
        red!75!black,
        line width=1.4pt
    },
    recovery/.style={
        red!75!black,
        line width=0.9pt
    },
    immune/.style={
        blue!75!black,
        line width=2.0pt,
        densely dotted
    }
]


\draw[site] (-1,0) -- (-1,6.4);
\draw[site] ( 0,0) -- ( 0,6.4);
\draw[site] ( 1,0) -- ( 1,6.4);
\draw[site] ( 2,0) -- ( 2,6.4);

\draw[-{Latex[length=2mm]},black]
    (-1.55,0) -- (-1.55,7.35);

\node[rotate=90,font=\small]
    at (-1.82,3.2) {time};

\node[below,font=\small] at (-1,0) {$x-1$};
\node[below,font=\small] at ( 0,0) {$x$};
\node[below,font=\small] at ( 1,0) {$x+1$};
\node[below,font=\small] at ( 2,0) {$x+2$};


\draw[infected]
    (-1,0.20) -- (-1,6.4);

\draw[infection]
    (-1,0.85) -- (0,0.85);

\node[above,font=\scriptsize]
    at (-0.5,0.85)
    {$e_0\in\mathcal R_m$};


\draw[infected]
    (0,0.85) -- (0,2.30);

\draw[infection]
    (0,1.65) -- (1,1.65);

\node[above,font=\scriptsize]
    at (0.5,1.65)
    {$e_1\in\mathcal A_m$};


\draw[infected]
    (1,1.65) -- (1,4.65);

\draw[infection]
    (1,2.45) -- (2,2.45);

\node[above,font=\scriptsize]
    at (1.5,2.45)
    {$e_2\in\mathcal A_m$};

\draw[infected]
    (2,2.45) -- (2,6.4);


\draw[recovery]
    (-0.10,2.20) -- (0.10,2.40);
\draw[recovery]
    (-0.10,2.40) -- (0.10,2.20);

\draw[immune]
    (0,2.40) -- (0,3.45);

\filldraw[
    orange!80!black,
    fill=white,
    line width=0.7pt
]
    (0,3.45) circle (1.7pt);


\draw[infection]
    (1,4.00) -- (0,4.00);

\node[above,font=\scriptsize]
    at (0.5,4.00)
    {$e_1^{\leftarrow}\in\mathcal R_{m+1}$};

\draw[infected]
    (0,4.00) -- (0,6.4);


\draw[recovery]
    (0.90,4.55) -- (1.10,4.75);
\draw[recovery]
    (0.90,4.75) -- (1.10,4.55);

\draw[immune]
    (1,4.75) -- (1,5.55);

\filldraw[
    orange!80!black,
    fill=white,
    line width=0.7pt
]
    (1,5.55) circle (1.7pt);


\draw[
    decorate,
    decoration={brace, amplitude=3pt}
]
    (-0.10,1.68) -- (-0.10,2.18);

\node[
    font=\scriptsize,
    anchor=east
]
    at (-0.17,1.93) {$R_x$};

\draw[
    decorate,
    decoration={brace, amplitude=3pt}
]
    (-0.10,2.42) -- (-0.10,3.43);

\node[
    font=\scriptsize,
    anchor=east
]
    at (-0.17,2.925) {$I_x$};

\draw[
    decorate,
    decoration={brace, amplitude=3pt}
]
    (-0.10,3.47) -- (-0.10,3.98);

\node[
    font=\scriptsize,
    anchor=east
]
    at (-0.17,3.725) {$A_{x+1,x}$};

\draw[
    decorate,
    decoration={brace, amplitude=3pt, mirror}
]
    (1.10,1.68) -- (1.10,4.53);

\node[
    font=\scriptsize,
    anchor=west
]
    at (1.17,3.105) {$R_{x+1}$};
\node[
    draw=black,
    fill=white,
    rounded corners=2pt,
    inner xsep=7pt,
    inner ysep=4pt,
    align=center,
    font=\scriptsize
]
    at (-0.3,5.78)
    {$R_x+I_x+A_{x+1,x}<R_{x+1}$\\[2pt]
     $\Longrightarrow$ reverse gate succeeds};

\end{tikzpicture}}
		\caption{Return-gate exploration. The return infection
$e_0\in\mathcal R_m$ generates first-infection edges
$e_1,e_2\in\mathcal A_m$. The figure illustrates
a possible return through the reverse gate associated
with $e_1$.}
		\label{fig:return-gate-exploration}
	\end{figure}
	
	Consequently, by independence and the exponential law,
the conditional probability that a candidate gate
succeeds is at most
\[
\begin{aligned}
q
&:=\mathbf P_{\lambda,\mu}
  (R_u+I_u+A_{v,u}<R_v)\\
&=\mathbf E_{\lambda,\mu}
  \left[e^{-(R_u+I_u+A_{v,u})}\right]\\
&=\frac12\,\frac{\mu}{1+\mu}\,
  \frac{\lambda}{1+\lambda}.
\end{aligned}
\]
The estimate is made at the creation of the gate,
before any information about its descendants is revealed. By symmetry, this bound applies to both orientations
of each candidate gate.
	
Charging each candidate gate at its creation time
and using the tower property,
\eqref{eq:number-gates} and
\eqref{eq:first-infection-count} yield
\[
\begin{aligned}
\mathbf E_{\lambda,\mu}Z_{m+1}
&\leq q\,\mathbf E_{\lambda,\mu}G_m\\
&\leq 2q\,\mathbf E_{\lambda,\mu}
\bigl[Z_m+|\mathcal A_m|\bigr]\\
&\leq
\frac{\mu\lambda(1+2\lambda)}
{(1+\mu)(1+\lambda)}
\mathbf E_{\lambda,\mu}Z_m.
\end{aligned}
\]
Writing
$m_{\mathrm{ret}}
=\mu\lambda(1+2\lambda)/
[(1+\mu)(1+\lambda)]$,
condition~\eqref{eq:extinction-criterion}
gives $m_{\mathrm{ret}}<1$.

Since the initial infected set is finite, the same
geometric exploration, started from its sites, gives
$\mathbf E_{\lambda,\mu}Z_1<\infty$. Consequently,
\[
\mathbf P_{\lambda,\mu}(Z_m>0)
\leq \mathbf E_{\lambda,\mu}Z_m
\leq m_{\mathrm{ret}}^{m-1}
\mathbf E_{\lambda,\mu}Z_1
\longrightarrow 0.
\]
Since $\mathbf E_{\lambda,\mu}Z_m<\infty$ for every
$m$, each generation is almost surely finite.
Consequently, infinitely many return infections would
require nonempty generations of arbitrarily large
order. But
\[
\mathbf P_{\lambda,\mu}(Z_m>0)
\longrightarrow 0,
\]
so there are almost surely only finitely many return
infections. By~\cite[Theorem~5]{Heeszel2026},
survival from a finite infected set is therefore
impossible, completing the proof.
\end{proof}

\appendix
\section{Moment estimates for regeneration cycles}
\label{app:moment-estimates}

We establish exponential moment bounds for the duration,
displacement, and maximal oscillation of regeneration
cycles. The latter controls deviations of the right edge
from its regenerative skeleton, allowing limit theorems
for the regeneration increments to extend to the full
right-edge process.

\begin{proof}[Proof of Proposition~\ref{prop:cycle-exponential-moments}]
We first consider a generic post-first cycle started from $Q$, using the
notation of Subsection~\ref{subsec:one-cycle-from-Q}. Thus $N$ is the first accepted trial and
\[
        \tau_+=T_N+s_*,
        \qquad
        X_+=a_N+1,
        \qquad
        M_+=\sup_{0\le t\le\tau_+}|r_t|.
\]
By Proposition~\ref{prop:barrier-regeneration-increments},
$(X_+,\tau_+,M_+)$ has the common law of
$(\Delta X_n,\Delta\tau_n,M_n)$ for $n\ge2$.

We first control the cost of rejected trials. Write
\[
        F_i:=\{\text{the $i$-th trial is rejected}\}.
\]
Conditionally on the information $\mathcal G_i$ available before the
$i$-th trial, Lemma~\ref{lem:freshness-adaptive-trials} gives the
fresh-trial law. A preparation failure has a uniformly bounded cost. If the
preparation succeeds and the launched toy process dies, we use the standard
restart consequence of
Proposition~\ref{prop:peierls-after-old-tail-deletion}: there exist $c,C>0$
such that
\begin{equation}
\label{eq:toy-late-extinction}
         \mathbf P_\mu\bigl(t<\tau^{\rm toy}<\infty\bigr)
        \le Ce^{-ct},
        \qquad t\ge0.
\end{equation}
Indeed, after each finite failed block cluster one restarts the block
construction from the surviving toy process, after the fixed local
seed-preparation step. The uniform positive survival probability and the
exponential tail for the height of a finite cluster give
\eqref{eq:toy-late-extinction} by the usual geometric restart argument.

Moreover, the toy process can be coupled in the obvious way with a classical contact process, and the CPTI right edge is bounded above by the right edge of a classical contact process under the common Harris construction. The standard
finite-speed estimate for the classical contact process therefore controls
both the spatial range revealed while the toy failure is resolved and any
simultaneous rightward CPTI displacement. Consequently there exist $c,C>0$,
uniformly in $i$, such that
\begin{equation}
\label{eq:rejection-cost-tail}
\mathbf P_Q\left(
        O_i\ge t,\ F_i\mid\mathcal G_i
\right)
\le Ce^{-ct},
\qquad t\ge0.
\end{equation}
Since a successful preparation is followed by toy survival with conditional
probability $\widehat\theta>0$, the fresh-trial construction also gives
\begin{equation}
\label{eq:moment-proof-uniform-success}
 \mathbf P_Q(F_i\mid\mathcal G_i)
\le 1-p_{\rm br}
\qquad\text{a.s.}
\end{equation}

We now control the accepted level. By
\eqref{eq:one-cycle-target-recursion},
\begin{equation}
\label{eq:bound-target-level}
        a_N\le N+\sum_{i=1}^{N-1}O_i.
\end{equation}
Choose $\alpha>0$ so small that \eqref{eq:rejection-cost-tail} implies
\[
K:=
\sup_i\left\|
 \mathbf E_Q\left[
        e^{\alpha O_i} \Ind{F_i}\mid\mathcal G_i
\right]
\right\|_\infty<\infty.
\]
For $0<\delta<\alpha$, conditional H\"older and
\eqref{eq:moment-proof-uniform-success} give
\[
\mathbf E_Q\left[
        e^{\delta(O_i+1)} \Ind{F_i}\mid\mathcal G_i
\right]
\le
e^\delta K^{\delta/\alpha}
(1-p_{\rm br})^{1-\delta/\alpha}.
\]
The right-hand side tends to $1-p_{\rm br}<1$ as
$\delta\downarrow0$. Hence, for $\delta>0$ sufficiently small, there is
$\chi<1$ such that
\begin{equation}
\label{eq:failed-trial-contraction}
 \mathbf E_Q\left[
        e^{\delta(O_i+1)} \Ind{F_i}\mid\mathcal G_i
\right]
\le \chi
\qquad\text{a.s.}
\end{equation}
On $\{N=n\}$, \eqref{eq:bound-target-level} yields
\[
        e^{\delta a_N}
        \le
        e^\delta\prod_{i=1}^{n-1}
        e^{\delta(O_i+1)} \Ind{F_i}.
\]
Iterating \eqref{eq:failed-trial-contraction} and summing over $n\ge1$
gives
\begin{equation}
\label{eq:record-exponential-moment}
         \mathbf E_Q[e^{\delta a_N}]<\infty.
\end{equation}

We next control the cycle duration. Write
\[
        r_t^{\rm toy}:=\sup\eta_t^{\{0\}},
        \qquad
        \bar r_t^{\rm toy}:=\sup_{0\le s\le t}r_s^{\rm toy},
\]
with the convention that $r_t^{\rm toy}=-\infty$ when
$\eta_t^{\{0\}}=\varnothing$. The same restart construction, together
with the linear growth of the infinite good-block path in
Proposition~\ref{prop:peierls-after-old-tail-deletion}, gives constants
$v_0,c,C>0$ such that
\begin{equation}
\label{eq:toy-lower-speed-on-survival}
\mathbf P_\mu\left(
        \tau^{\rm toy}=\infty,\,
        \bar r_t^{\rm toy}\le v_0t-C
\right)
\le Ce^{-ct},
\qquad t\ge0.
\end{equation}
Under $Q$ the distinguished toy process is conditioned on survival. Since
every toy infection path is a legal CPTI infection path,
\eqref{eq:toy-lower-speed-on-survival} therefore implies, after changing
constants,
\begin{equation}
\label{eq:Q-running-record-lower-tail}
         \mathbf P_Q(\bar r_t\le v_0t-C)\le Ce^{-ct},
        \qquad t\ge0.
\end{equation}
By construction, $a_N>\bar r_{\rho_{N-1}}$. Hence $a_N$ has not been
reached before the $N$-th trial and
\[
        T_N=\inf\{t\ge0:r_t=a_N\}.
\]
Since the right edge moves by nearest-neighbor jumps,
\begin{equation}
\label{eq:record-passage-inclusion}
\{T_N>t\}
\subseteq
\{a_N>v_0t-C\}
\cup
\{\bar r_t\le v_0t-C\}.
\end{equation}
Combining \eqref{eq:record-exponential-moment},
\eqref{eq:Q-running-record-lower-tail}, and
\eqref{eq:record-passage-inclusion}, we obtain
\begin{equation}
\label{eq:record-passage-tail}
          \mathbf P_Q(T_N>t)\le Ce^{-ct},
        \qquad t\ge0.
\end{equation}
Since $\tau_+=T_N+s_*$ and $X_+=a_N+1$, after decreasing $\delta$ if
necessary,
\begin{equation}
\label{eq:generic-time-space-exponential-moments}
         \mathbf E_Q[e^{\delta\tau_+}]
        + \mathbf E_Q[e^{\delta|X_+|}]<\infty.
\end{equation}

It remains to control the oscillation during the cycle. Choose
$\varepsilon>0$ small. For the rightward excursion, the CPTI is dominated
by the usual contact process started from the infected half-line, so the
standard finite-speed estimate gives an exponential bound. For the
leftward excursion, if $r_s\le-m$, then the distinguished barrier, which is
infected at time $s$, must also lie at or to the left of $-m$. Since this
barrier is a toy infection path and hence a classical graphical infection
path, the same finite-speed estimate applies. Conditioning on toy survival
only changes the multiplicative constant. Thus there exist $c,C>0$ such
that
\begin{equation}
\label{eq:Q-short-time-excursion-tail}
 \mathbf P_Q\left(
        \sup_{0\le s\le\varepsilon m}|r_s|\ge m
\right)
\le Ce^{-cm},
\qquad m\ge0,
\end{equation}

Moreover,
\begin{equation}
\label{eq:cycle-oscillation-inclusion}
\{M_+\ge m\}
\subseteq
\{\tau_+>\varepsilon m\}
\cup
\left\{\sup_{0\le s\le\varepsilon m}|r_s|\ge m\right\}.
\end{equation}
By \eqref{eq:record-passage-tail},
\eqref{eq:Q-short-time-excursion-tail}, and
\eqref{eq:cycle-oscillation-inclusion},
\[
         \mathbf P_Q(M_+\ge m)\le Ce^{-cm},
        \qquad m\ge0,
\]
and hence, after decreasing $\delta$ once more,
\begin{equation}
\label{eq:generic-oscillation-exponential-moment}
         \mathbf E_Q[e^{\delta M_+}]<\infty.
\end{equation}
Equations \eqref{eq:generic-time-space-exponential-moments} and
\eqref{eq:generic-oscillation-exponential-moment}, together with
Proposition~\ref{prop:barrier-regeneration-increments}, prove
\eqref{eq:cycle-exponential-moments}.

For the first cycle, the rejected-trial estimate
\eqref{eq:rejection-cost-tail} and the contraction
argument remain valid. The block construction from
the initially infected half-line and the classical
finite-speed comparison give
\begin{equation}
\label{eq:first-cycle-inputs}
\mathbf P_\mu(\bar r_t\le v_0t-C)\le Ce^{-ct},
\qquad
\mathbf P_\mu\left(
\sup_{0\le s\le\varepsilon m}|r_s|\ge m
\right)\le Ce^{-cm}.
\end{equation}
Writing $\tau_1=T^{(1)}+s_*$ and $X_1=a^{(1)}+1$,
the preceding arguments, with
\eqref{eq:first-cycle-inputs} replacing the
corresponding $Q$-estimates, yield
\[
\mathbf E_\mu[e^{\delta\tau_1}]
+\mathbf E_\mu[e^{\delta|X_1|}]
+\mathbf E_\mu[e^{\delta M_1}]<\infty,
\]
proving \eqref{eq:first-cycle-moments}.
\end{proof}

\section*{Acknowledgements}

R.S. is supported by CNPq grant 403423/2023-6. C.T. is supported by CNPq grant 151432/2024-4. M.E.V. is partially supported by CNPq grant 310734/2021-5 and by FAPERJ grant E-261200.473/2026. 

\section*{Statement on AI use}
Artificial intelligence (AI) tools were used solely for linguistic proofreading and assistance with the code used to produce figures. The authors reviewed and verified all AI-assisted changes and code. AI tools were not employed in the formulation of mathematical results or the development of proofs. The authors take full responsibility for the mathematical content and the final manuscript.

\bibliographystyle{plainurl}
\bibliography{refbib}
\end{document}